\documentclass[11pt, leqno, twoside]{amsart}
\usepackage{latexsym}
\usepackage{amssymb}
\usepackage{amsmath}
\usepackage{bbm}
\usepackage{bm}

\usepackage{tikz}
\usepackage{tkz-euclide}
\usetikzlibrary{decorations.fractals}
\usetikzlibrary{decorations.footprints}
\usepackage[colorlinks=true, pdfstartview=FitV, linkcolor=red, citecolor=red, urlcolor=blue, backref=page]{hyperref} 

\usepackage{tikz,amsthm,amsmath,amstext,amssymb,amscd,epsfig,euscript,pspicture,multicol,graphpap,graphics,graphicx,enumerate,subfig,sidecap,wrapfig,color,pict2e}
\usepackage{mathrsfs} 
\usepackage{upgreek}
\usepackage[utf8]{inputenc}

\usepackage[top=1in, bottom=1in, left=1in, right=1in, a4paper, marginparwidth=60pt]{geometry}

\usepackage{marginnote}
\calclayout
\allowdisplaybreaks

\theoremstyle{plain}
\newtheorem{theorem}[equation]{Theorem}
\newtheorem{lemma}[equation]{Lemma}
\newtheorem{corollary}[equation]{Corollary}
\newtheorem{proposition}[equation]{Proposition}

\theoremstyle{definition}
\newtheorem{definition}[equation]{Definition}

\theoremstyle{remark}
\newtheorem{remark}[equation]{Remark}

\numberwithin{equation}{section}

\renewcommand*{\backref}[1]{}
\renewcommand*{\backrefalt}[4]{%
 \ifcase #1 (Not cited.)%
   \or        (Cited on page~#2.)%
    \else      (Cited on pages~#2.)%
    \fi}

\DeclareMathOperator{\diam}{diam}
\DeclareMathOperator{\dv}{div}
\DeclareMathOperator{\Tr}{Tr}

\DeclareMathOperator{\dist}{dist}

\DeclareMathOperator{\conv}{conv}
\DeclareMathOperator{\ext}{ext}

\DeclareMathAlphabet{\mathsfit}{T1}{\sfdefault}{\mddefault}{\sldefault}
\SetMathAlphabet{\mathsfit}{bold}{T1}{\sfdefault}{\bfdefault}{\sldefault}

\DeclareRobustCommand{\SkipTocEntry}[5]{}

\newcommand{\set}[2]{\left\{#1 : #2\right\}}
\newcommand{\emp}{\emptyset}
\newcommand{\sub}{\subseteq}
\newcommand{\mns}{\setminus}
\newcommand{\N}{\mathbb{N}}
\newcommand{\R}{\mathbb{R}}

\newcommand{\br}[1]{\overline{#1}}
\newcommand{\del}{\partial}

\newcommand{\I}{\mathbb{I}}
\newcommand{\eps}{\varepsilon}
\newcommand{\inv}[1]{{#1}^{-1}}
\newcommand{\dx}{\, dx}
\newcommand{\loc}{\text{\rm loc}}
\newcommand{\Om}{\Omega}
\newcommand{\om}{\omega}
\newcommand{\inp}[2]{\big\langle #1,#2\big\rangle}
\newcommand{\vertiii}[1]{{\left\vert\kern-0.25ex\left\vert\kern-0.25ex\left\vert #1 
    \right\vert\kern-0.25ex\right\vert\kern-0.25ex\right\vert}}
\newcommand{\gr}{\nabla}
\newcommand{\wto}{\rightharpoonup}
\newcommand{\lap}{\Delta}
\newcommand{\h}{\mathcal{H}}
\newcommand{\on}{%
  \,\raisebox{-.127ex}{\reflectbox{\rotatebox[origin=br]{-90}{$\lnot$}}}\,}

\newcommand{\Sn}{\mathbb S}
\newcommand{\Snn}{\Sn^{n-1}}
\newcommand{\g}{\mathbf{g}}
\newcommand{\F}{\mathscr{F}}
\newcommand{\LL}{\mathscr{L}}
\newcommand{\cof}{\mathcal{C}}
\newcommand{\K}{\mathcal{K}}
\newcommand{\kr}{\upkappa}
\newcommand{\cov}{\text{\raisebox{2pt}{$\bigtriangledown$}}}
\newcommand{\covtwo}{{\cov}{}^2}
\usepackage{filecontents}

\begin{document} 

	\allowdisplaybreaks[2]

\title[Minkowski Problem for $p$-harmonic measures]{On the Minkowski Problem for $p$-harmonic measures}

\author[M. Akman]{Murat Akman}
\address{{Murat Akman}\\
Department of Mathematical Sciences
\\
University of Essex
\\
Colchester CO4 3SQ, United Kingdom}	\email{murat.akman@essex.ac.uk}

\author[S. Mukherjee]{Shirsho Mukherjee}
\address{{Shirsho Mukherjee}\\
Department of Mathematical Sciences
\\
University of Essex
\\
Colchester CO4 3SQ, United Kingdom}	\email{shirsho.mukherjee@essex.ac.uk, \ m.shirsho@gmail.com}

\thanks{Both authors have been supported by the EPSRC New Investigator award [grant number EP/W001586/1].}

\date{\today}


\begin{abstract}
 
We study the Minkowski problem corresponding to the $p$-harmonic measures and obtain results previously known for 
harmonic measures due to Jerison \cite{J1}. We show that a class of Borel measures on spheres can be prescribed by $p$-harmonic measures 
on convex domains. 
 \end{abstract}

\maketitle
\setcounter{tocdepth}{2}
\tableofcontents


\section{Introduction}

One of the most fundamental problems in the geometric analysis of convex bodies is the Minkowski problem. The original form of the problem
seeks to find a convex polyhedron that has a specified set of vectors and real numbers as its normals and 
surface areas of its faces under stipulated conditions. It can be reformulated more generally as seeking to find a convex domain that induces a specified measure on a sphere as the image of its surface measure mapped by its normals. In this form, the original problem is identical to the case of having 
a discrete measure on a sphere. Similar related problems are to show the uniqueness of such domains up to translations and to find convex domains with smooth boundaries when the regularity of the Gaussian curvature of its boundary as a function of its normal is specified. 

The literature goes back to Minkowski \cite{M2,M1} who 
proved the existence and uniqueness of a convex polyhedron with specified discrete spherical measure if and only if the centroid is at the origin and support of the measure is not concentrated on an equator (or great sub-spheres). The problem for measures having continuous density was also illustrated by Minkowski and later solved for general measures by Alexandrov \cite{A1, A2} and Fenchel-Jessen \cite{FJ}. 
The problem corresponding to convex domains with smooth boundaries involves establishing regularity of solutions for Monge-Amp\`ere type equations on the sphere, which is rather difficult and has involved many significant works including 
Lewy \cite{Lewy}, Pogorelov \cite{Pog}, Nirenberg \cite{Nir}, 
Cheng-Yau \cite{CY}, Caffarelli \cite{Caff90-2,Caff90-1,Caff91}, etc. Furthermore, the harmonic measures on the boundary of convex domains being equivalent to the surface measures (see Dahlberg \cite{D}), 
an analogous Minkowski problem with surface measures being replaced by harmonic measures, has been studied by Jerison \cite{J1}. 
The objective of the present paper is to study the Minkowski problem corresponding to $p$-harmonic measures for $1<p<\infty$, which, to our knowledge, remained open.

For $n\geq 3$, we consider convex domains $\Om\subset \R^n$ and the Gauss map $\g_\Om :\del\Om\to \Snn$ where $\g_\Om(x)$ is the outer unit normal at $x\in \del\Om$, which is defined almost everywhere 
on the boundary as 
convexity implies that $\del\Om$ is locally Lipschitz. We consider a positive finite Borel measure $ \mu $ on $ \mathbb{S}^{n-1}$ 
that satisfies the following necessary conditions: 
\begin{equation}\label{eq:excond} 
\begin{aligned}
(i)&\, \,   { \int_{ \mathbb{S}^{n-1}} } | \inp{\zeta}{ \xi} | \, d \mu (  \xi )  
>  0, \quad\forall\   \zeta \in \mathbb{S}^{n-1} ,\\  
(ii)& {  \int_{ \mathbb{S}^{n-1}} }  \xi \, d \mu ( \xi )  = 0. \\
\end{aligned}
\end{equation}
The original Minkowski problem seeks the existence up to a translation of a unique compact convex domain $\Om$ with non-empty interior such that $(\g_\Om)_* \h^{n-1}\on_{\del\Om} = \mu$, where 
$(\g_\Om)_* \h^{n-1}$ is the pushforward measure on $\Snn$ and $\mu$ satisfies 
$(i)$ and $(ii)$ of the above. We show the 
existence of the domain corresponding to $\h^{n-1}$ being replaced by $p$-harmonic measures. 

Given any convex set $K\subset \R^n$, the $p$-harmonic measure $\om_p$ supported on $\del K$ associated to a function $u=u_K$ is given by
 $$\om_p(E) = \int_{\del K\cap E} |\gr u|^{p-1}d\h^{n-1},$$ 
for any measurable $E\sub \R^n$, where $u$ vanish on $\del K$ and is non-negative and $p$-harmonic in a neighborhood of it. 
This includes $K$ being open or the closure of an open convex domain as well as the case of $K$ being a convex set of empty interior. 
It also includes the case of $u$ having a possible blow-up or being arbitrarily large somewhere in the interior of the domain away from the boundary, 
that resembles the behavior of Green's function for the linear case. 
The $p$-harmonic measure has been studied by Lewis \cite{Lnote}, Lewis-Nystr\"{o}m-Vogel \cite{LNV}, see also 
\cite{ALV, HKM, LLN, LN}, etc. 
The Gauss map $\g_K$, which is defined almost everywhere from $\del K$ to $\Snn$ when $K$ is convex and of non-empty interior with $\g_K(x)\in \Snn$ being the outer unit normal for $\h^{n-1}$-a.e. $x\in \del K$, can be generalized to be defined as a set-valued map where $\g_K(x)$ constitutes the set of all normals 
at any $x\in \del K$ as in \cite{J1}; the notion of push forward $(\g_K)_*$ for measures can be accordingly generalized, which coincides with the case when $\g_K$ is single-valued. This is necessary to include convex sets of empty interior and lower dimensional sets. Further details are in Section \ref{sec:prelim}. 

The main result of this paper is the following. 
\begin{theorem}\label{thm:main}
Given a finite regular Borel measure $ \mu $ on $\mathbb{S}^{n-1}$ satisfying conditions \eqref{eq:excond}, there exists a bounded convex domain $\Om$ with non-empty interior such that $(\g_\Om)_* \om_p\on_{\del\Om} = \mu$, where $\g_\Om$ is the $\h^{n-1}$-a.e. defined Gauss map and 
$\om_p$ is any $p$-harmonic measure on $\del\Om$ for $1<p<\infty$. 
\end{theorem}

Many further problems of similar character as the Minkowski problem have been studied with different measures as replacement of the surface measures, e.g. capacitary measures, curvature measures, etc. We refer to \cite{BKLYZ,CW,GG,GLL,Guan-Ma,L-O}, and references therein.  
The proof of the existence of domains involves variational methods with respect to the Minkowski sum of convex domains $\Om+t\Om'$ as in Jerison \cite{J1, J2} and Colesanti et al. \cite{CNSXYZ} (see also \cite{ALSV,AGHLV,ALVcalvarpde}); 
adaptations of such methods form the key elements of our proof of existence. The problem of the existence of domains is converted to a constrained minimization problem of functions by virtue of the dualistic relation between convex domains to their support functions and 
positive continuous functions to their Wulff shapes, see Section \ref{sec:prelim} for details. 
However, contrary to surface measures, harmonic measures, 
capacitary measures, and others, the $p$-harmonic measures are not uniquely defined on the boundary of convex domains as their definition in general involves the choice of a $p$-harmonic function. Nevertheless, we show that the arbitrarity of the choice can be maintained as the characteristics of being $p$-harmonic at the boundary are similar for all $p$-harmonic measures and they admit natural approximants that can be uniquely defined for the purpose of establishing existence via limiting argument. 

The paper is organized as follows. In Section \ref{sec:prelim}, we develop our notations and 
provide a self-contained account of preliminaries on convex domains, geometric features of their boundaries, and approximations. We also enlist 
the definitions and known results on $p$-harmonic functions, $p$-harmonic measures, and Minkowski problems. In Section \ref{sec:ex}, 
we compute the first variation of densities of $p$-harmonic measures defined on smooth convex domains of the form $\Om+t\Om'$. 
Then we provide an account of weak convergence of $p$-harmonic measures which are used to construct approximants of the measures. 
The properties of the first variation are used in the variational methods for the existence of the domains for discrete measures, which are further used as approximations for the existence of the domains for general measures in the proof of Theorem \ref{thm:main} at the end. 

The uniqueness up to translation and regularity of the domains are nontrivial and have not been addressed in this paper. 
These objectives are to be pursued in the future. 




\section{Notations and Preliminaries}\label{sec:prelim}

In this section, we introduce some notations and review some well-known properties of convex domains,  $p$-harmonic measures and Minkowski problems. 

For $n\geq 1$, we denote points in $\R^n$ as 
$ x = ( x_1,
 \dots,  x_n) $ so that if $\{e_1,\ldots,e_n\}$ is the standard basis then $( x_1,
 \dots,  x_n)=\sum_{i=1}^n x_ie_i$. The standard inner
product on $ \mathbb R^{n} $ shall be denoted by $\inp{\cdot}{\cdot} $ and the Euclidean norm by $|\cdot|$. 
For functions $f:\R^n\to \R$ and $F= (f_1, \ldots, f_m):\R^n\to \R^m$, the gradient defined by $\gr f = \sum_{i=1}^n (\del_{i} f) e_i$, the Jacobian 
defined by $ DF = \sum_{i=1}^m\sum_{j=1}^n (\del_{i} f_j) e_i\otimes e_j$ and 
the Hessian defined by
$ D^2f =D(\gr f)= \sum_{i,j=1}^n (\del^2_{i,j} f) e_i\otimes e_j$, are as usual.  

In general, hyperplanes in $\R^n$ shall be denoted as 
\begin{equation}\label{eq:hyp}
H_{ y,\alpha} := \set{x\in \R^n}{\inp{x}{ y}= \alpha},
\end{equation}
 for some $ y\in \R^n$ and $\alpha\in \R$; the half-spaces corresponding to $\inp{x}{ y}< \alpha$ and $\inp{x}{ y}> \alpha$ shall be respectively denoted as 
 $H_{ y,\alpha}^-$ and $H_{ y,\alpha}^+$. Thus, we have $H_{ y,\alpha}^-=H_{ -y,-\alpha}^+,\ H_{ y,\alpha}^+=H_{-y,-\alpha}^-$ and 
 $H_{ y,\alpha}=H_{-y,-\alpha}$, for any $ y\in \R^n$ and $\alpha\in \R$. 
 In particular, $H_{e_n, 0}^+$ is denoted by $\R^n_+$. 

Given a measure space $(X, \mu)$ and 
a map $f:X\to Y$, the pushforward measure $f_*\mu$ on $Y$ is defined on any measurable subset $ E\sub Y$ as
$$ (f_* \mu)(E) = \mu(\inv{f}(E)),$$
which is absolutely continuous with respect to $\mu$ and in the infinitesimal form, for any function $g:Y\to \R$ it is written as 
$g\, d(f_* \mu)= (g\circ f) d\mu$. 
For any Borel set $E\sub \R^n$, the 
$k$-dimensional Hausdorff measure on $\R^n$ is defined by 
\[
\h^{k}(E)=\lim_{\delta\to 0^+} \inf\Big\{\sum_{j} r_j^{k} \, : \, E\subset\bigcup_j B_{r_j}(x_j), \, \, r_j\leq \delta\Big\},
\] 
where $B_r(x)= \{y\in \R^n : |y-x|<r\}$ is the standard metric ball of radius $r$ centered at $x\in \R^n$. 
The distance of $E$ from a point $y \in \R^n$ is defined by 
$\dist(y, E)=\inf\{|x-y|\, :\, x\in E\}$. 
The boundary of the unit sphere $B_1(0)$ of $\R^n$ shall be denoted as $ \Snn$. 
We have the transformation to polar coordinates  
$\R^n\mns\{0\} \to (0,\infty)\times \Snn$ as $x\mapsto (|x|, x/|x|)$ and its inverse
$(r, \theta) \mapsto r \theta $. 
For any $ F\in L^1(\R^n) $ we have 
\begin{equation}\label{eq:polarint}
\int_{\R^n} F(x) dx = \om_n \int_0^\infty \bigg[\int_{\Snn} F(r \theta)r^{n-1}d \theta\bigg] dr,
\end{equation}
where the Lebesgue measure $\mathcal L^n$ of $\R^n$ and the uniform measure $\h^{n-1}$ on $\Snn$
are abbreviated as $dx$ and $d \theta$ in infinitesimal form and $\om_n=2\pi^{n/2}/\Gamma(n/2)$ is the surface area of $\Snn$. 

\subsection{Convex Domains}\label{subsec:convdom}
Here we shall denote $\Om\subset \R^n$ as a convex domain, i.e. open and connected non-empty convex subset and $K\subset \R^n$ as a convex body, i.e. closure of a bounded convex domain, hence a compact convex subset with non-empty interior. We highlight some essential notions and properties related to convex sets and functions. 
For further details on the theory of convex bodies, we refer to Schneider \cite{Sc} and  Gardner \cite{G}. 

For $A\subset \R^n$, a function $f: A \to \R$ is convex iff the epigraph $K= \{(x,c)\in A\times \R: f(x)\leq c\}$ is convex,  
which implies $f$ is continuous and locally Lipschitz; we denote $\K(A)$ as the set of all convex functions $f: A \to \R$ and hence, $\K(A)\subset C^{0,1}_\loc (A)$. 
We recall that $f\in C^1(\Om)\cap \K(\Om)$ iff 
$\inp{\gr f(x)-\gr f(y)}{x-y}\geq 0$ for all $x,y\in\Om$ and $f\in C^2(\Om)\cap \K(\Om)$ iff $D^2f \geq 0$. In general, for $f\in \K(\R^n)$, the \textit{sub-differential} of $f$ at $x\in \R^n$, given by 
\begin{equation}\label{eq:subdiff}
\del f(x)= \{ v\in \R^n : f(y)\geq f(x) + \inp{v}{y-x}\},
\end{equation}
is a convex subset and $f$ is differentiable at $x$ iff $\del f(x)$ is a singleton. An important class of convex functions, is the class of \textit{sub-linear functions} that are \textit{sub-additive} and 
\textit{homogeneous} of degree $1$; in other words, they are defined by functions $h:\R^n\to \R$ satisfying 
\begin{equation}\label{eq:sublinear}
h(x+y) \leq h(x) +h(y) \quad \text{and}\quad h(\lambda x) = \lambda h(x)\qquad \forall\ x,y\in \R^n,\ \lambda\geq 0. 
\end{equation}
A closed subset $K\subset \R^n$ is convex iff for any $z\in \R^n$ there exists a unique projection $p_K(z)\in K$ satisfying 
$|z-p_K(z)| = \dist(z, K)$; in this case, we have $\dist(\cdot, K) \in \K(\R^n)\cap C^1(\R^n\mns K)$ and 
\begin{equation}\label{eq:gradist}
 \xi_K(z):=\gr\dist(z, K) = \frac{z-p_K(z)}{\dist(z, K)},\qquad\forall\ z\in \R^n\mns K.
\end{equation}
Thus, for any $z\in\R^n\mns K$, we have $p_K(z)\in \del K$ and $|\gr\dist(z, K)|=1$. For any $x\in K$ and $t\in [0,1]$, since $t\mapsto |z-(tx+(1-t)p_K(z))|^2$ attains minimum at $t=0$, 
its slope at $t=0$ is non-negative which leads to 
\begin{equation}\label{eq:projcond}
\inp{x-p_K(z)}{z-p_K(z)} \leq 0, \qquad\forall\ x\in K,\ z\in \R^n\mns K.
\end{equation}
Furthermore, the ray emerging from $p_K(z)$ connecting $z$ being $R_{z,K}=\{p_K(z)+\lambda  \xi_K(z): \lambda\geq 0\}$, we note that $p_K(z')=p_K(z)$ and $ \xi_K(z')= \xi_K(z)$ for all $z'\in R_{z,K}$. The continuity of $p_K$ also ensures that for any $x\in \del K$ their exists $z\in \R^n\mns K$ such that $p_K(z)=x$ and hence $R_{z,K}\sub \inv{p_K}(x)$. 
If $\del K$ is smooth enough at $x\in \del K$ then $R_{z,K}= \inv{p_K}(x)$ and $ \xi_K(z)$ coincides with the outer unit normal at $x\in \del K$.  
This leads to the following. 
\begin{definition}\label{def:normal}
If $K$ is closed and convex, for any $x\in \del K$, the vector $ \xi_K(z)\in \Snn$ as in \eqref{eq:gradist} is called a (outer) \textit{normal vector} of $K$ at $x$ if $z\in \inv{p_K}(x)$.
\end{definition} 
For any $A\subset \R^n$ , the \textit{support function} $h_A: \R^n \to \R$ is defined as 
\begin{equation}\label{eq:suppf}
h_A( y) = \sup \set{\inp{x}{ y} }{ x\in A}.
\end{equation}
Thus, we have $h_\Om = h_{\bar\Om}$ and $h_K( y) = \max_{x\in K}\inp{x}{ y}$ for a compact convex subset $K\subset\R^n$. Note that $h_K$ is linear (i.e. $h_K( y)= \inp{x}{ y}$ for every 
$ y\in \R^n$) iff $K=\{x\}$. 
Furthermore, the convexity implies $x\in K$ iff $\inp{x}{ y}\leq h_K( y)$ for every $y$, i.e. 
in terms of the subdifferential \eqref{eq:subdiff}, $\del h_K(0)= K$. The \textit{supporting hyperplane} of $K$ with outer normal $ y\in  \R^n\mns\{0\}$, is given by 
\begin{equation}\label{eq:supphp}
H_K ( y):= H_{ y,h_K( y)} = \set{x\in \R^n}{\inp{x}{ y}= h_K( y)}, 
\end{equation}
so that $K \subset H_{ y,h_K( y)}^-=\set{x\in \R^n}{\inp{x}{ y}\leq h_K( y)}$ (called \textit{supporting half-space}) and 
it can be shown in terms of the subdifferential \eqref{eq:subdiff} that for any $y\in \R^n\mns\{0\}$ we have $$\del h_K( y)= K\cap H_K( y)\subset \del K,$$ see \cite{Sc}. 
Thus, $h_K$ is differentiable at $ y\in \R^n\mns\{0\}$ with $\gr h_K( y)=x$ iff $K\cap H_K( y)=\{x\}$; in this case $h_{K\cap H_K( y)}(\cdot)= \inp{x}{\cdot}= \inp{\gr h_K( y)}{\cdot}$ and $K$ is called \textit{strictly convex} at $x$ and if it holds for every point then the set is called strictly convex. For such domains, since $h_K$ is homogeneous of degree $1$, hence $ y\mapsto \gr h_K( y)$ is homogeneous of degree $0$, i.e. 
\begin{equation}\label{eq:gradhom}
\gr h_K (\lambda  y) = \gr h_K( y), \qquad \forall\ \lambda\geq 0.
\end{equation}
In general, it is evident that support functions are sub-linear and moreover, all sub-linear functions are support functions by virtue of the following theorem, see \cite[Theorem 1.7.1]{Sc}. 
\begin{theorem}\label{thm:subsupp}
For any sublinear function $h:\R^n\to \R$, there exists a unique convex body $K=\{x\in \R^n : \inp{x}{ y}\leq  h( y)\ \forall\  y\in\R^n\}$ such that $h_K= h$. 
\end{theorem}
Due to homogeneity, it suffices to define the set $K=\{x\in \R^n : \inp{x}{ \xi}\leq  h( \xi)\ \forall\  \xi\in\Snn\}$ in Theorem \ref{thm:subsupp} having $h$ as its support function and henceforth, it suffices to consider support functions restricted to $\Snn$ for all things and purposes. Also, this implies that for any convex set $K$, we have
\begin{equation}\label{eq:kinthk}
K = \bigcap_{ \xi\in \Snn} \{x\in \R^n : \inp{x}{ \xi}\leq  h_K( \xi)\}. 
\end{equation}
Thus, it is evident that we have $K\sub K'$ iff $h_K\leq h_{K'}$. 
Now, alongside having supporting hyperplanes with a fixed normal, one can have a supporting hyperplane at a fixed point in the boundary. 
These two notions are related to each other via projection.
Indeed, 
from \eqref{eq:projcond}, we note that $h_K(z-p_K(z))=\inp{p_K(z)}{z-p_K(z)}$ for any $z\in \R^n\mns K$ and hence, we have  
\begin{equation}\label{eq:supproj}
h_K( \xi_K(z))=\inp{p_K(z)}{ \xi_K(z)}, 
\end{equation}
with $ \xi_K(z)$ is as in \eqref{eq:gradist}. 
Therefore, for any $x\in \del K$, their exists a supporting hyperplane of $K$ at $x$ given by 
\begin{equation}\label{eq:supphp1}
H_K[x] = \set{x'\in \R^n}{\inp{x'-p_K(z)}{z-p_K(z)}= 0}, \quad\text{with}\ z\in \inv{p_K}(x);
\end{equation}
which is related to \eqref{eq:supphp} by $H_K[x]=H_K(z-p_K(z))$ from \eqref{eq:supproj}. Thus, from homogeneity of support functions, $H_K[x]= H_K( \xi_K(z))$ 
for any $z\in \inv{p_K}(x)$ where $ \xi_K(z)$ is as in \eqref{eq:gradist} and moreover, $\del K$ is smooth at $x\in \del K$ iff 
$H_K[x]$ is unique with a unique outer unit normal.  

We recall that, for any $E,F\sub \R^n$, the \textit{Minkowski sum} is defined by 
\begin{equation}\label{eq:minksum}
 E  +  F  = \{ x + y :  x \in E, y \in F   \}; 
\end{equation}
also, we shall denote by $ c E =  \{ c y : y \in E\} $ for any $c\in \R$. It is not hard to check that if $E, F$ are convex then so is $E+F$ and $\alpha E$ for all $\alpha\geq 0$. Furthermore, for compact convex sets, the decomposition to Minkowski sums is unique, i.e. if $E+K_1=E+K_2$ then $K_1=K_2$. 
In particular, note that $B_r(x)= \{x\}+ r B_1(0)$. One of the most important features of convex sets is that their support functions are Minkowski additive, i.e. for compact convex sets $E, F\subset \R^n$, we have 
\begin{equation}\label{eq:suppminkadd}
h_{\alpha E+\beta F} = \alpha h_E + \beta h_F \qquad\forall\ \alpha, \beta\geq 0, 
\end{equation}
and hence, for supporting hyperplanes as in \eqref{eq:supphp}, $H_{\alpha E+\beta F}(\cdot) = \alpha H_E (\cdot)+ \beta H_F(\cdot)$ holds as well. 
For any general $E\subset \R^n$, the \textit{convex hull} of $E$, denoted by $\conv(E)$, is defined as the smallest convex set (intersection of all convex sets) containing $E$ and we have $$\conv(E+F)= \conv(E)+\conv(F).$$ 
If $E$ is convex, then $\conv(E)=E$, if $E$ open (or closed) so is $\conv(E)$. In fact, points in convex hull of a set comprise convex combinations of affinely independent points of the set; i.e. for any $x\in \conv(E)$, there exists affinely independent $x_1,\ldots, x_k\in E$ for $k\leq n+1$ and $\lambda_1,\ldots, \lambda_k\geq 0$ such that $x=\sum_j \lambda_j x_j$ where $\sum_j \lambda_j=1$ 
(see \cite[Theorem 1.1.4]{Sc}). 
\begin{definition}\label{def:convhull}
The convex hulls of finitely many points are called \textit{polytopes}. For any $m\leq n$, if the points $x_1,\ldots, x_{m+1}\in \R^n$ are affinely indepedent then the polytope 
$P= \conv(\{x_1,\ldots, x_{m+1}\})$ is called a $m$-\textit{simplex} with vertices $x_1,\ldots, x_{m+1}$.
\end{definition} 
Thus, $\conv(E)$ is the union of all $m$-simplices for $m\leq n$ with vertices in $E$. Conversely, for a convex set $K$, the set of points in $ K$ that cannot be written as $\lambda x+ (1-\lambda)y$ for any $x,y\in K$ is denoted as $\ext(K)$ with elements called \textit{extreme points}. Any compact convex set is the closed convex hull of its extreme points, i.e. $K= \br{\conv(\ext(K))}$ (Krein-Milman theorem) and hence, polytopes are the only convex bodies with finitely many extreme points. Therefore, for any polytopes $P_1$ and $P_2$ and $\alpha,\beta\geq 0$, evidently $\alpha P_1+\beta P_2$ is also a polytope. 
If $H$ is a supporting hyperplane of a polytope $P=\conv(\{x_1,\ldots,x_k\})$, then we have 
$$ H\cap P = \conv(H\cap \{x_1,\ldots,x_k\}),$$
which is also a polytope of a lower dimension; such polytopes are called \textit{faces} of $P$. In fact, all polytopes are finite intersections of their supporting 
half-spaces corresponding to their faces codimension $1$, see \cite[Theroem 2.4.3]{Sc}. Furthermore, we have the following characterization theorem of polytopes, we refer to \cite[Corollary 2.4.4 and Theorem 2.4.6]{Sc} for the proof. 
\begin{theorem}\label{thm:charpoly}
A convex set is a bounded intersection of finitely many closed half-spaces iff it is a polytope. Moreover, there exists $ \xi_1, \ldots,  \xi_k\in \Snn$ with respect to which a polytope $P$ is uniquely represented as $P = \big\{x\in \R^n : \inp{x}{ \xi_i}\leq h_P( \xi_i)\ \forall\ i\in\{1,\ldots, k\}\big\}$. 
\end{theorem}
The vectors $ \xi_1, \ldots,  \xi_k\in \Snn$ of Theorem \ref{thm:charpoly} are normal vectors of $P$ with $h_P( \xi_i)$'s uniquely determining the polytope and are chosen so that 
$P\cap H_P( \xi_i)$ are the $1$-codimensional faces. 
\subsection{Convex Geometry}\label{subsec:convgeom} 
We enlist some essential features and geometric notions of convex domains with smooth boundaries which have some generalizations to general convex domains. 

Recalling the notion of smoothness of boundaries, $\Om\subset \R^n$ is said to be of class $C^k$ (resp. $C^{k,\alpha}$) for $k\in \N$ and $0<\alpha\leq 1$ if for any 
$x\in \del \Om$, there exists $r>0$ and a bijection $\psi: B_r(x)\to \R^n$ with $\psi, \inv{\psi} \in C^k$ (resp. $C^{k,\alpha}$) such that 
$\psi (B_r(x)\cap \Om) \subset \R^n_+$ and $\psi (B_r(x)\cap \del\Om) \subset \del\R^n_+$; equivalently, for any $x\in \del \Om$ there exists a neighborhood $U$ and $\phi\in C^{k}(\R^{n-1})$ (resp. $C^{k,\alpha}(\R^{n-1})$) such that 
after a possible rotation, we have
\begin{equation}\label{eq:omlocal}
\begin{aligned}
\Om\cap U &= U\cap \set{(x_0,x')\in\R\times \R^{n-1}}{\phi(x')<x_0};\\
\del\Om\cap U &= U\cap \set{(x_0,x')\in\R\times \R^{n-1}}{\phi(x')=x_0}.
\end{aligned}
\end{equation}
In addition, if $\Om$ is convex then $\phi$ of \eqref{eq:omlocal} is a convex function. If $\Om$ is of class $C^k$ (resp. $C^{k,\alpha}$), then any 
$f\in C^k(\del\Om)$ (resp. $C^{k,\alpha}(\del\Om)$) can be extended to a function in $C^k(\bar\Om)$ (resp. $C^{k,\alpha}(\bar\Om)$). 

If $\del K$ is of class $C^1$ for a convex domain $K$, then for every $x\in \del K$ there exists a normal vector at $x$ as in Definition \ref{def:normal} that is unique up to scaling. This gives rise to  
$\g_K : \del K\to \Snn$, called \textit{Gauss map}, where $\g_K(x)$ is the outer unit normal at $x\in \del K$. Thus, notice that for any $z\in \inv{p_K}(x)$, we have 
$\g_K(x) = (z-x)/|z-x|= \xi_K(z)$ with $ \xi_K(z)$ is as in \eqref{eq:gradist} and in terms of support function, from \eqref{eq:supproj}, we have
\begin{equation}\label{eq:supg}
h_K (\g_K(x)) = \inp{x}{\g_K(x)}, \quad \forall\ x\in\del K. 
\end{equation}
In particular, for balls we have $\g_{B_r(x_0)} : \del B_r(x_0)\to \Snn$ given by $\g_{B_r(x_0)}(x) = (x-x_0)/r$ for any $x_0\in \R^n$ and $r>0$.  
Recall that the tangent space at $x\in\del K$ is given by 
$$ T_x(\del K)= \set{y\in \R^n}{\inp{y}{\g_K(x)}=0} $$
and from \eqref{eq:supg}, we observe that $H_K(\g_K(x))= \{x\}+T_x(\del K)$. The boundary $\del K$ being a $C^1$-manifold, we have $\g_K\in C^1(\del K, \Snn)$ and we recall some familiar notions from differential geometry. 
Since $T_x(\del K)=T_{\g_K(x)}(\Snn)$, the following \textit{Weingarten map} is defined by 
$$\mathcal{W}_K(x):= d\,\g_K(x) : T_x(\del K)\to T_x(\del K).$$ The \textit{principal curvatures} are the eigenvalues 
$\kappa_1(x),\ldots, \kappa_{n-1}(x)$ of $\mathcal{W}_K(x)$; the \textit{mean curvature} and the \textit{Gaussian curvature} are respectively given by
$$ H(x)=\frac{\Tr(\mathcal{W}_K(x))}{n-1}
= \frac{1}{n-1}\big(\kappa_1+\ldots+ \kappa_{n-1}\big)\quad \text{and}\quad \kr(x)=\det(\mathcal{W}_K(x))= \kappa_1\ldots\kappa_{n-1}.$$  
If $K$ is strictly convex then $K\cap H_K(\g_K(x))= \{x\}$ for any $x\in \del K$ and $\inv{\g_K}:\Snn\to \del K$ is well defined; moreover, 
$h_K$ being differentiable everywhere in $\R^n\mns\{0\}$, from \eqref{eq:supg} we have 
\begin{equation}\label{eq:hgrh}
h_K( \xi) = \inp{ \xi}{\inv{\g_K}( \xi)} \quad\text{and}\quad \gr h_K( \xi) = \inv{\g_K}( \xi), \quad\forall\  \xi\in \Snn. 
\end{equation}

\begin{definition}\label{def:c+}
A convex set $K$ is said to be of class $C^k_+$ (resp. $C^{k,\alpha}_+$) for any $k\in \N,  \alpha\in (0,1]$ if $\del K$ is of class $C^k$ (resp. $C^{k,\alpha}$) and the Gauss map 
$\g_K : \del K\to \Snn$ is a diffeomorphism. 
\end{definition}
The Gauss map being a diffeomorphism is equivalent to the Weingarten map being of maximal rank everywhere. This is further equivalent to all principal curvatures being non-zero or equivalently the Gaussian curvature being non-zero; the convexity implies that the curvatures are positive everywhere. The corresponding domains are called \textit{strongly convex}. 

If $K$ is of class $C^2_+$, note that $\inv{\g_K}\in C^1(\Snn, \del K)$ and $h_K\in C^2(\R^n\mns\{0\})$. Hence, 
from homogeneity of $h_K$, \eqref{eq:hgrh} and \eqref{eq:gradhom}, we have 
\begin{equation}\label{eq:hgradh}
h_K( y) = \inp{ y}{\gr h_K( y)}\quad \text{and}\quad D^2h_K( y)  y =0, \qquad\forall\  y\in \R^n\mns\{0\},
\end{equation}
where the latter equality of the above can be obtained by differentiating the former. We have the inverse Weingarten map 
$$\inv{\mathcal{W}_K} ( \xi):= \inv{(\mathcal{W}_K(\inv{\g_K}( \xi)))}= d\,\inv{\g_K}( \xi) : T_ \xi (\Snn) \to T_ \xi (\Snn)$$ defined for all $ \xi\in \Snn$, which is the non-singular part of $D^2h_K( \xi)$. 
Indeed, with $ \xi\in \Snn$ being fixed, there exists an orthonormal basis $\{e^1,\ldots,e^{n-1},  \xi\}$ of $\R^n$ and hence, $\{e^i\}$ span the tangent space 
$T_\xi (\Snn)$. Any $x\in \R^n$, in this basis, is given by 
\begin{equation}\label{eq:basis}
x= \sum_{k=1}^{n-1}x^k e^k + \inp{x}{ \xi} \xi,\qquad\text{with}\quad x^k= \inp{x}{e^k}.
\end{equation}
Since $D^2h_K ( \xi)  \xi =0$ from \eqref{eq:hgradh}, the only
non-zero entries of $D^2h_K ( \xi)$ are $\inp{D^2 h_K( \xi) e^j}{e^i}$ for $i,j\in\{1, \ldots, n-1\}$ which, from \eqref{eq:hgrh}, are 
the entries of $d \gr h_K( \xi)= d\,\inv{\g_K}( \xi)=\inv{\mathcal W_K}( \xi)$, with respect to the above basis. In other words, we have 
\begin{equation}\label{eq:d2wing}
D^2 h_K( \xi) = \sum_{i,j=1}^{n-1} \inp{\inv{\mathcal W_K}( \xi) e^j}{e^i} e^i \otimes e^j. 
\end{equation}
Thus, $D^2 h_K( \xi)$ has eigenvalues $\{1/\kappa_1(\inv{\g_K}( \xi)),\ldots, 1/\kappa_{n-1}(\inv{\g_K}( \xi)), 0\}$, where $\kappa_i$'s are the principal curvatures of $\del K$. With $\xi\in U\subset \Snn$ and a coordinate chart $\varphi:U \to V\subset \R^{n-1}$, the covariant derivatives of $f: \Snn \to \R$ of first and second orders are locally defined by 
\begin{equation}\label{eq:cov}
\begin{aligned}
(i)&\ \cov f:=\sum_{i=1}^{n-1}(\cov_i f) \, e^i, \quad \text{where}\quad \cov_if (x):= \del_i (f\circ \inv{\varphi})(\varphi(x)),\\
(ii)&\ \covtwo f:=\sum_{i,j=1}^{n-1} (\cov_{i,j} f)\, e^i \otimes e^j,  \quad \text{where}\quad \cov_{i,j}f (x):= \del_{i,j} (f\circ \inv{\varphi})(\varphi(x)),
\end{aligned}
\end{equation}
Letting $\chi=\inv{\varphi}: V\to U\subset\Snn$, note that since $|\chi|^2=1$ hence by differentiating successively we have 
$\inp{\chi}{\del_j\chi}=0$ and $\inp{\chi}{\del_{i,j}\chi}=-\inp{\del_i\chi}{\del_j\chi}$. We can choose $U, \varphi$ such that
$\chi(0)=\xi$ and that $\del_j\chi(0)= e^j$ for all $j\in \{1,\ldots, n-1\}$ taking 
$U=\Snn_+=\big\{x\in \R^n: \inp{x}{\xi}=\sqrt{\scriptstyle{1-\sum_{i=1}^{n-1}|x^i|^2}} \big\}$ and $\varphi=(x^1,\ldots, x^{n-1})$, leading to $\chi (z_1,\ldots, z_{n-1})= 
\sum_{i=1}^{n-1}z_i e^i +\sqrt{\scriptstyle{1-\sum_{i=1}^{n-1}|z_i|^2}} \, \xi$. 
Thus, we have $\del_j\chi = e^j- (z_j/ \sqrt{\scriptstyle{1-\sum_{i=1}^{n-1}|z_i|^2}} )\, \xi$ and $\del_{i,j}\chi(0)=-\delta_{i,j}\chi(0)$  
(one can define a similar chart on the other hemisphere for $U=\Snn_-$ corresponding to the negative square root). 
Taking $h_K$ restricted to $\Snn$, \eqref{eq:hgradh} can be written in terms of local coordinates as 
\begin{equation}\label{eq:hgradhcov}
h_K\circ \chi =\inp{\chi}{\gr h_K\circ \chi}\quad\text{and}\quad (D^2 h_K\circ \chi)\chi =0. 
\end{equation}
By differentiating and using the above, we get
$$\del_j(h_K\circ \chi)
= \inp{\del_j\chi}{\gr h_K\circ \chi}+\inp{\chi}{(D^2 h_K\circ \chi)\del_j\chi}=\inp{\del_j\chi}{\gr h_K\circ\chi}.$$
Evaluating at $0$, we get $\cov_j h_K(\xi) = \inp{\gr h_K(\xi)}{e^j}$. Hence, from \eqref{eq:basis} and \eqref{eq:hgradhcov}, notice that 
\begin{equation}\label{eq:grhcov}
\gr h_K ( \xi)= h_K( \xi) \xi + \cov h_K( \xi). 
\end{equation}
By differentiating twice at $0$ and using the choice of chart, we get 
\begin{align*}
\del_{i,j}(h_K\circ \chi)(0)&= \inp{\del_{i,j}\chi(0)}{\gr h_K\circ \chi(0)}+ \inp{\del_j\chi(0)}{(D^2 h_K\circ \chi)(0) \del_i\chi(0)}\\
&= -\delta_{i,j} \inp{\chi(0)}{\gr h_K(\chi(0))}+ \inp{e^j}{(D^2 h_K(\chi(0)) e^i}
\end{align*}
which, from \eqref{eq:hgradh} and \eqref{eq:d2wing}, further implies 
\begin{equation}\label{eq:cov2h}
\covtwo  h_K ( \xi) = \inv{\mathcal W_K}( \xi) - h_K( \xi)\I,
\end{equation}
where $\I$ is the identity matrix. 
Therefore, we have the following for $C^2_+$ domains; 
\begin{equation}\label{eq:detkinv}
 \det\big(\covtwo  h_K( \xi) + h_K( \xi) \I\big) = \det (\inv{\mathcal W_K}( \xi))= 1/\kr(\inv{\g_K}( \xi)), 
\end{equation}
where $\kr$ is the Gaussian curvature. One can show that a domain is $C^2_+$ iff it is $C^2$ and $\kr >0$. 
From transformation rule of the Jacobian, note that
 $$(\g_K)_*\h^{n-1}\on_{\del K}= |\det(d\,\inv{\g_K})| \h^{n-1}\on_{\Snn}=1/(\kr\circ \inv{\g_K}) \h^{n-1}\on_{\Snn}$$ and from \eqref{eq:detkinv}, for any $f \in L^1(\del K, \h^{n-1})$, we have 
\begin{equation}\label{eq:intbd}
\int_{\del K} f(x) \, d\h^{n-1}(x)= \int_{\Snn} f(\inv{\g_K}( \xi)) \det\big(\covtwo  h_K( \xi) + h_K( \xi) \I\big) \, d \xi. 
\end{equation}

In general, 
since convex functions are locally Lipschitz, the boundary of a convex domain $\Om$ can be locally written as a graph of a Lipschitz function and the Gauss map can be defined $\h^{n-1}$-a.e. on 
$ \del \Om$. Precisely, for a small enough neighborhood $U$, if $\phi\in \K(\R^{n-1})\cap C^{0,1}(\R^{n-1})$ is as in \eqref{eq:omlocal}
then, for a.e. $x'\in \R^{n-1}$, $\phi$ is differentiable and 
$$g_\Om(\phi(x'), x') = (-1,\gr \phi(x'))/\sqrt{1+|\gr \phi(x')|^2}.$$ Moreover, a classical theorem of Aleksandrov \cite{BF, A3} (see also \cite{EG}) is the following. 
\begin{theorem}[Aleksandrov]\label{thm:alex}
If $\phi:\R^n\to \R$ is convex, then $\phi$ is twice-differentiable a.e.
\end{theorem}
Consequently, for a convex domain $\Om$ and $K=\bar\Om$, we have $\mathcal W_K(x)$ and $\kr(x)\geq 0$ well-defined for $\h^{n-1}$-a.e. $x\in \del \Om$ and
the surface measure $S_\Om$ on $\Snn$ induced by the almost everywhere defined Gauss map 
is defined by 
\begin{equation}\label{eq:surfm}
S_\Om(E) := \h^{n-1}\Big(\set{x\in \del \Om}{ \g_\Om(x)\ \text{is well-defined},\ \g_\Om(x)\in E} \Big)
\end{equation}
for any measurable $E\sub \Snn$ and we denote $S_\Om=S_K$ if $K=\bar\Om$. For any $w\in L^1(\Snn, \h^{n-1})$, the integral formula 
\eqref{eq:intbd} generalizes to 
\begin{equation}\label{eq:intbdgen}
\int_{\del \Om} w(\g_\Om(x)) \, d\h^{n-1}(x)= \int_{\Snn} w( \xi) \, d S_\Om ( \xi).
\end{equation}
It is evident that if $g_\Om$ is defined everywhere then $S_\Om=(\g_{\Om})_*\h^{n-1}$ and if $\Om$ is $C^2_+$, then $d S_\Om ( \xi)= \det (\covtwo  h_\Om( \xi) + h_\Om( \xi) \I)  d \xi$ as earlier.

For a general convex set $K$ the Gauss map is re-defined as a set-valued map at $x\in \del K$ as 
\begin{equation}\label{eq:gmap}
\g_K(x)=\set{\xi\in \Snn}{\inp{x'-x}{\xi}\leq 0, \ \forall\ x'\in K},
\end{equation}
i.e. the set of all normal vectors $\xi\in \Snn$ of every supporting hyperplane $H_K[x]$ at $x\in \del K$.  
The push forward $(\g_K)_*$ can be re-defined with respect to the re-defined inverse image 
\begin{equation}\label{eq:ginv}
\inv{\g_K}(E)=\set{x\in \del K}{\g_K(x)\cap E\neq \emp},
\end{equation}
so that we can re-define
$S_K=(\g_{K})_*\h^{n-1}$ as $S_K(E)=\h^{n-1}(\inv{\g_K}(E))$ for any convex domain $K$ and measurable $E\sub \Snn$. 
If $K=\bar\Om$ for a convex domain $\Om$ then, the set \eqref{eq:gmap} is a singleton for $\h^{n-1}$-a.e. $x\in \del K$.
More generally, the set \eqref{eq:gmap} can be a cone at non-smooth boundary points or a discrete set if $K$ is of lower dimension. 
Henceforth, according to the regularity of the convex set $K$, the notions $\g_K$ and $S_K$ shall be assumed at the appropriate level of generality, hereafter. 

For any $A\subset \R^n$, let us denote $C_+(A)$ as the class of positive, continuous functions $f:A\to \R$. If $0\in \Om$ then from \eqref{eq:kinthk} 
it is clear that $h_\Om \in  C_+(\Snn)$. Conversely, any $h\in C_+(\Snn)$ coincides with the support function of a domain $\Om_h$ except possibly for a set of zero surface measure on $\Snn$; 
the domain $\Om_h$, called the \textit{Wulff shape} of $h$, is defined by 
\begin{equation}\label{eq:wshape}
\bar\Om_h = \bigcap_{ \xi\in \Snn} \{x\in \R^n : \inp{x}{ \xi}\leq  h( \xi)\}.
\end{equation}
Indeed, since $h$ is positive and continuous, it is evident that $\Om_h$ is a bounded convex domain with $0\in \Om_h$. One can extend $h$ to be a 1-homogeneous function on $\R^n\mns\{0\}$ as $ y\mapsto | y| h( y/| y|)$ and letting $K=\bar\Om_h$, clearly $h_K\leq h$. Moreover, if $x\in \del K$ then 
\eqref{eq:wshape} implies $\inp{x}{ \xi}= h( \xi)$ for some $ \xi\in \Snn$ and \eqref{eq:kinthk} implies $\inp{x}{ \xi'}= h_K( \xi')$ for some $ \xi'\in \Snn$; if $h( \xi)>h_K( \xi)$ then $ \xi\neq  \xi'$ and hence $H_{ \xi, h( \xi)}$ and $H_{ \xi', h_K( \xi')}$ are two different supporting hyperplanes at $x$ with normals $ \xi,  \xi'$ making $x$ a non-smooth point of $\del K$. Thus Aleksandrov's theorem enforces $h_K( \xi)=h( \xi)$ for  
$S_K$-a.e. $ \xi\in \Snn$ and thus we have a generalization of Theorem \ref{thm:subsupp} as follows. 
\begin{theorem}\label{thm:possupp}
Given any $h\in C_+(\Snn)$, there exists a bounded convex domain $\Om=\Om_h$ as in \eqref{eq:wshape} such that $0\in \Om$ and 
$h_\Om ( \xi)= h( \xi)$ holds for $S_\Om$-a.e. $ \xi\in \Snn$. 
\end{theorem}

Let $\Om$ be a compact convex domain with $0\in \Om$ and $K=\bar \Om$. Then, integral formulae on $\del K$ can be obtained independently using the \textit{radial map}
$r_K: \R^n\mns \{0\}\to [0,\infty)$ defined by 
\begin{equation}\label{eq:radialmap}
r_K( y)= \sup\set{r\geq 0}{r y\in K}. 
\end{equation}
It is clear that for any $\alpha>0$, we have $r_{\alpha K}( y)= \alpha r_K( y)$ and $r_K(\alpha y)= r_K( y)/\alpha$. Also, it is not hard to see that  
$r_K( y) y\in \del K$ for any $ y \in \R^n\mns \{0\}$. 
This gives rise to the \textit{radial projection} $\rho_K: \Snn\to \del K$ defined by $\rho_K( \theta)=r_K( \theta) \theta$, which is the polar coordinate image 
on $\del K$. From the polar coordinate formula \eqref{eq:polarint} we have the following, for any $F\in L^1(K)$ 
\begin{equation}\label{eq:polark}
\int_{K} F(x) dx = \om_n \int_{\Snn} \bigg[\int_0^{r_K( \theta)}F(r \theta)r^{n-1}dr\bigg] d \theta;
\end{equation}
in particular, we have $|K|= (\om_n/n)\int_{\Snn} r_K( \theta)^n d \theta$. Hence, from the divergence theorem, we have that 
$\int_{\del\Om}\inp{x}{\g_K(x)} d\h^{n-1} = \int_{\Snn} |\rho_K( \theta)|^n d \theta$, and 
$$\h^{n-1}\on_{\del K} = (\rho_K)_*\Big(|x|^n/\inp{x}{\g_K(x)} \, \theta\on_{\Snn}\Big)$$
so that, $d\h^{n-1} = |\rho_K( \theta)|^n/\inp{\rho_K( \theta)}{\g_K(\rho_K( \theta))} d\theta = |r_K( \theta)|^{n-1}/\inp{\theta}{\g_K(\rho_K( \theta))} d\theta$, 
in infinitesimal form. Thus we have 
the following due to Aleksandrov 
(see \cite[Lemma 1]{A2}). 
\begin{theorem}\label{thm:intbdal}
Let $\Om$ be a bounded convex domain with $0\in \Om$, $K=\bar\Om$ and $\g_K$ be the Gauss map defined $\h^{n-1}$-a.e. on $\del K$. 
For any $f \in L^1(\del K, \h^{n-1})$,  we have the following; 
\begin{equation}\label{eq:intbdal}
\int_{\del K} f(x) \, d\h^{n-1}(x)= \int_{\Snn} f\big(\rho_K( \theta)\big)\frac{r_K( \theta)^{n-1}}{\inp{\g_K(\rho_K( \theta))}{ \theta}}  d \theta. 
\end{equation}
\end{theorem}
Furthermore, as $\h^{n-1}$-a.e. points on $\del K$ are smooth, hence recalling \eqref{eq:supg}, we can also re-write the above formulae using 
 $d\h^{n-1} = |\rho_K( \theta)|^n/h_K (\g_K(\rho_K( \theta))) \, d\theta$. Unlike \eqref{eq:detkinv},\eqref{eq:intbd} above, this density can be used for any 
 convex body that are not necessarily smooth. 
 
\subsection{Hausdorff distance and approximations}\label{subsec:hausd}
Here we provide some details on properties of distance functions on convex sets, results on compactness and approximations. 

The \textit{Hausdorff distance} between Borel sets $ E, E' \subset\mathbb R^n $ 
is defined as 
 \[ 
  d_{\h} ( E, E' ) = \max  \Big(  \sup   \{ \dist ( y, E ) : y \in E' \}\, ,\,   \sup \{ \dist ( y, E' ) : y \in E \} \Big); \]  
 equivalently, we have $ d_{\h} ( E, E' ) = \sup_{y\in \R^n} |\dist ( y, E )-\dist ( y, E' )|$. 
 The Hausdorff distance of convex sets can be characterized by support functions as 
\begin{equation}\label{eq:hdh}
d_\h (E, F)=d_\h (\del E, \del F) = \|h_E -h_F\|_{L^\infty(\Snn)},
\end{equation}
see \cite[Lemma 1.8.1, Lemma 1.8.14]{Sc}. 
Any compact convex set $K$ can be approximated by polytopes with respect to $d_\h$, since, given any $\eps>0$, the compactness implies there exists points $x_1,\ldots, x_N$ such that $K$ contained in the union of $B_\eps (x_j)$'s and then the polytope $P = \conv(\{x_1,\ldots, x_{N}\})$ satisfies $P\subset K\subset P+ \eps \,\Snn$ and $d_\h(K, P)<\eps$. More generally, the space of bounded convex bodies is locally compact with respect to the Hausdorff distance.
\begin{theorem}[Blaschke selection theorem]\label{thm:bla}
Given any bounded sequence of convex sets $\{K_m\}$, there exists a subsequence $\{K_{m_j}\}$ and a convex set $K$ such that $d_\h (K_{m_j}, K)\to 0^+$ as $j\to \infty$. 
\end{theorem}

 If $d_\h (K_j,K)\to 0^+$ as $j\to \infty$, \eqref{eq:hdh} and Aleksandrov's theorem implies $\g_{K_j}\to \g_K$ $\h^{n-1}$ a.e. Moreover, 
since $\dist(x_0,\del K) \leq \dist(x_0,\del K_j)+ d_\h (K_j,K)$, it is not hard to see from local compactness, that then any subsequential limit $x_0$ of radial projections $\rho_{K_j}(\theta)$ lies on $\del K$, leading to uniform convergence of the projections, hence $r_{K_j}\to r_K$ uniformly as $j\to \infty$ and the densities of 
$\h^{n-1}$ with respect to radial projections
\begin{equation}\label{eq:dradconv}
 \frac{|\rho_{K_j}( \theta)|^n}{h_{K_j} (\g_{K_j}(\rho_{K_j}( \theta)))} \to \frac{|\rho_K( \theta)|^n}{h_K (\g_K(\rho_K( \theta)))}\qquad \h^{n-1}-a.e.\ \theta\in \Snn.
\end{equation}

We have the following approximation theorem for convex domains, see Weil \cite{W1,W2}. 
\begin{theorem}\label{thm:c2approx}
Given any convex set $K$ with $\h^{n-1}$-a.e. defined Gaussian curvature $\kr$, there exists a sequence $\{K_j\}$ of convex sets
of class $C^2_+$ such that as $j\to \infty$ we have $d_\h(K_j,K)\to 0^+$ and if $\kr_j$ is the Gaussian curvature of $K_j$ then the following holds: 
\begin{enumerate}
\item $1/\kr_j\to 1/\kr$ pointwise $\h^{n-1}$-a.e. as $j\to \infty$;
\item $1/\kr_j\to 1/\kr$ in $L^1(\Snn, \h^{n-1})$ as $j\to \infty$;
\item $\int_{\cdot} (1/\kr_j) \, d \xi \wto \int_{\cdot} (1/\kr) \, d \xi $ as $j\to \infty$.
\end{enumerate}
\end{theorem}
Thus, from Theorem \ref{thm:c2approx}, we can conclude that given a convex domain $\Om$ with Gaussian curvature $\kr$, there exists $C^2_+$ (or $C^{2,\alpha}_+$) domains 
$\Om_j$ with Gaussian curvature $\kr_j$, such that as  $j\to \infty$, we have $d_\h(\Om_j,\Om)\to 0^+$ and $(\g_{\Om_j})_*\h^{n-1}\wto (\g_{\Om})_*\h^{n-1}$, i.e. 
for any $w\in C(\Snn)$, 
\begin{equation}\label{eq:wconvmsr}
\lim_{j\to \infty}\int_{\Snn} \frac{w( \xi)}{\kr_j(\inv{\g_{K_j}}( \xi))} d \xi= \int_{\Snn} \frac{w( \xi)}{\kr(\inv{\g_K}( \xi))} d \xi. 
\end{equation}
For general convex sets, we have the weak convergence of surface measures with respect to Hausdorff distance, i.e. as $j\to \infty$, 
\begin{equation}\label{eq:wcons}
S_{K_j} \wto S_K, \qquad\qquad \text{if}\quad d_\h(K_j,K)\to 0^+,
\end{equation}
when the measures are defined with respect to \eqref{eq:gmap} and \eqref{eq:ginv}, see \cite[Theorem 4.2.1]{Sc}. 
In general, there are results of weak convergence with respect to Hausdorff distance for any support measures and curvature measures for general 
convex sets, the proofs involve approximation of convex sets by convex domains, in particular, polytopes with non-empty interiors that have 
combinatorial formulae for support measures. We refer to \cite[Chaper 4,5]{Sc} for further details. 

\subsection{$p$-harmonic functions and measures}\label{subsec:pharm}  
Here we provide some preliminaries and some well-known classical results for $p$-harmonic functions and $p$-harmonic measures. 
We refer to \cite{HKM, LN, LLN, T, Lnote, LNV}, etc. for more details. 

Just as a harmonic function in a domain $\Om\subset \R^n$ being minimizer of the Dirichlet energy $\int_\Om |\gr w|^2 \dx $ with $w|_{\del\Om}=f$,  is the unique solution of the Dirichlet problem $ \lap v=0 $ in $\Om$ and $v=f$ on $\del\Om$, the minimizers of the $p$-Dirichlet energy 
$\int_\Om |\gr w|^p \dx $ are weak solutions to the $p$-Laplacian equation $\lap_p u= \dv(|\gr u|^{p-2}\gr u)=0$ in $\Om$, i.e. 
$$ \int_\Om |\gr u|^{p-2}\inp{\gr u}{\gr \phi}\dx = 0 ,\qquad \forall\ \phi\in C^\infty_0(\Om),$$
and are called $p$-\textit{harmonic functions} for $1\leq p<\infty$. They coincide with harmonic functions for $p=2$. The existence of weak solution $u\in W^{1,p}(\Om)$ is classical and follows from direct methods of calculus of variations. We have the following monotonicity inequality for any measurable $E\sub \R^n$, functions $u,v \in W^{1,p}(E)$ and a constant $c=c(p)>0$, given by 
\begin{equation}\label{eq:pmon}
\int_{E}\inp{ |\gr u|^{p-2} \gr u - |\gr v|^{p-2}\gr v}{ \gr u-\gr v}\dx\geq c
\begin{cases}
  {\scriptstyle \int_{E} |\gr u -\gr v|^p\dx} \ &\text{if}\ p\geq 2,\\
  \frac{\big(\int_{E} |\gr u -\gr v|^p\dx\big)^{2/p}}{\big(\int_{E}(|\gr u|+|\gr v|)^p\dx\big)^{2/p-1}}\ &\text{if}\ 1<p<2,
\end{cases}
\end{equation}
which can be used to show the comparison principle for the $p$-Laplacian, i.e. $\lap_p v\leq \lap_p u$ on $\Om$ and 
$u\leq v$ on $\del\Om$ in the trace sense implies $u\leq v$ on $\Om$; we refer to \cite{HKM} for a proof. 
The uniqueness of weak solutions of $ \lap_p u=0 $ in $\Om$ and $u=f$ on $\del\Om$ follows easily from the comparison principle. The regularity theory of $p$-harmonic functions is more involved. It has been shown by DiBenedetto \cite{Dib}, Lewis \cite{L2} and Tolksdorff \cite{T} that the weak solutions of $ \lap_p u=0 $ in $\Om$ for $1<p<\infty$ are locally $C^{1,\beta}$, i.e. there exists 
$\beta =\beta(n,p)\in (0,1)$ such that for any $B \subset\subset \Om$, we have
\begin{equation}\label{eq:c1reg}
\|u\|_{C^{1,\beta}(B)}= \|u\|_{L^\infty(B)} +\|\gr u\|_{L^\infty(B, \R^n)} 
+  \underset{\underset{x \neq y}{x,y\in B}}{\sup}  \frac{|\gr u(x)-\gr u(y)|}{|x-y|^\beta} \leq c\, \|u\|_{W^{1,p}(\Om)},
\end{equation}
for some $c=c(n,p,\diam(B))>0$. For $p>2$ the regularity is optimal. Furthermore, the continuity of the gradient implies that if $\gr u\neq 0$ in $\Om'\subset\Om$, then 
we can conclude $u\in C^\infty (\Om')$ 
from Schauder estimates, see \cite{GT}. The boundary regularity is also known for $\lap_p u= 0$ in $\Om$ and $u=f$ on $\del\Om$ if 
$\Om$ is of class $C^{1,\alpha}$ and $f\in C^{1,\alpha}(\del\Om)$, see Lieberman \cite{Li}; in this case $u\in C^{1,\beta}(\bar\Om)$ for some $\beta =\beta(n,p,\alpha)\in (0,1)$ along with the following global estimate
\begin{equation}\label{eq:c1regbd}
\|u\|_{C^{1,\beta}(\Om)} \leq C\big(n,p, \alpha, \|f\|_{C^{1,\alpha}(\del\Om)}, \|u\|_{W^{1,p}(\Om)}, \Om\big). 
\end{equation}

Given a bounded and sufficiently regular domain $\Om\subset \R^n$, for any $f \in C(\del\Om)$ and $v$ being the solution of the Dirichlet problem $ \lap v=0 $ in $\Om$ and $v=f$ on $\del\Om$, we have 
\begin{equation}\label{eq:harmsr}
v(x)= \int_{\del\Om} f(y) d\om^x(y) 
\end{equation}
from maximum principle and Riesz representation
theorem, 
where $\om^x$ is a measure on $\del\Om$ referred to as the \textit{harmonic measure} at $x\in\Om$. The measure can be prescribed by the Green's function with a pole at $x$, i.e. the weak (distributional) solution of the 
Dirichlet problem $ \lap G^x= \delta_x $ in $\Om$ and $G^x=0$ on $\del\Om$. If $\del\Om$ is smooth enough with unit normal $\nu$, it is not difficult to see that $d\om^x= 
(\del G^x/\del \nu) d\h^{n-1} =\inp{\gr G^x}{\nu} d\h^{n-1}$. 
This notion can be generalized to the case for the $p$-Laplacian for $1<p<\infty$. Given a neighbourhood $N$ of $\del\Om$, if $u\in W^{1,p}(\Om\cap N)$ is a positive weak solution to the
$p$-Laplacian in $\Om\cap N$, then upon zero extension $u\in W^{1,p}(N)$. Since $p$-superharmonic functions form a non-negative distribution as shown in \cite{HKM}, from a theorem of Schwartz, we have that 
there exists a non-negative Radon measure $\om_p$ on $\del\Om$ such that 
$$ \int_\Om |\gr u|^{p-2}\inp{\gr u}{\gr \phi}\dx = -\int_{\del\Om} \phi \, d\om_p ,$$
for any $\phi \in C^\infty_0(N)$; the measure $\om_p$ is called the $p$-\textit{harmonic measure} associated to $u$. 
Such measures can also be defined for more general $\mathcal A$-harmonic functions that are referred as Riesz measures, see \cite{HKM, KZ}.
If $\del\Om$ is Lipschitz with unit normal $\nu_{\del\Om}$ defined almost everywhere on $\del\Om$, then $d\om_p= -|\gr u|^{p-2}\inp{\gr u}{\nu_{\del\Om}} d\h^{n-1}$; moreover, since $ \{u> 0\}$ in $\Om$ 
we have $\del\Om= \del\{u>0\}$ and $\nu_{\del\Om}= -\gr u/|\gr u| $ and hence 
\begin{equation}\label{eq:phm}
d\om_p= |\gr u|^{p-1} d\h^{n-1}\on_{\del\Om}.
\end{equation}
Evidently, for the case $p=2$, if $N=\R^n\mns\{x\}$ and $u=G^x$ then $\om_2=\om^x$. More generally, the function $u$ may have a similar blow-up inside 
$\Om\mns N$ but the $p$-harmonicity in $\Om \cap N$ implies that $u\in C^{1,\beta}_\loc(\Om \cap N)$ for some $\beta\in(0, 1)$ from \eqref{eq:c1reg}. Given such a neighbourhood, it is evident that we can select finer reduced neighborhood $N'$ of $\del\Om$ i.e. $\del\Om\subset N'\sub N$, having the same properties. 
Such neighborhood $N$ (without relabelling) when chosen up to possible reduction, we can assume without loss of generality that 
$\gr u \neq 0$ in $\Om\cap N$ so that $u\in C^\infty (\Om \cap N)$ as stated above, and furthermore, we can assume 
$$\|u\|_{L^\infty(\del N\cap \Om)}+\|\gr u\|_{L^\infty(\del N\cap \Om)} <\infty ;$$
in other words, all possible singularities of $u$ are strictly in the interior of $\Om\mns \bar N$. Any such reductions can be chosen as long as it contains $\del\Om$ without changing any properties of the measure since the support of the measure is within $\del\Om$. 

The fundamental solution of the $p$-Laplacian, given by 
\begin{equation}\label{eq:fundp}
\Phi(x) =
\begin{cases}
1/|x|^\frac{n-p}{p-1} \quad &\text{for}\ p\neq n;\\
\log(1/|x|) \quad &\text{for}\ p=n,
\end{cases}
\end{equation}
solves $\lap_p \Phi=0$ in $\R^n\mns\{0\}$ and can be used to form examples of $p$-harmonic measures. 

The above notion of $p$-harmonic measure is defined for any open connected domains $\Om$, including convex domains. 
The notion extends to a convex set $K$ of non-empty interior in the same way taking $\Om$ as the interior of $K$ so that 
given a neighbourhood $N$ of $\del K$ and a function $u\in W^{1,p}(K\cap N)$, we have the notion of a $p$-harmonic measure 
$\om_p$ associated to $u$ given by $d\om_p= |\gr u|^{p-1} d\h^{n-1}\on_{\del K}$. 
Then, this can be further extended naturally to the 
case of $K$ being a convex set of empty interior i.e. $K=\del K$; given $N$ is a neighborhood of $K$ and 
$u\in W^{1,p}(N)$ being a positive weak solution to the
$p$-Laplacian in $N\mns K$ and $u$ vanish on $K$, there exists a non-negative Radon measure $\om_p$ on $K$ such that 
$$ \int_N |\gr u|^{p-2}\inp{\gr u}{\gr \phi}\dx = -\int_{K} \phi \, d\om_p ,$$
for any $\phi \in C^\infty_0(N)$ similarly as above which is called the $p$-harmonic measure associated to $u$, 
and $d\om_p= |\gr u|^{p-1} d\h^{n-1}\on_{ K}$ if $K$ is Lipschitz. We refer to 
\cite{LNld} for further details on quasi-linear equations on low dimensional sets. 
In general, since for all convex sets $K$ we have 
$\del K$ locally Lipschitz, hence $d\om_p= |\gr u|^{p-1} d\h^{n-1}\on_{\del K}$ for $u=u_K$.

\subsection{Minkowski problem}\label{subsec:minkprob}
Given a positive finite Borel measure $ \mu $ on $ \mathbb{S}^{n-1}$ that 
satisfies the conditions \eqref{eq:excond}, 
the Minkowski problem can be regarded as the inquiry for existence of a convex domain $\Om$ satisfying 
\begin{equation}\label{eq:minkprob}
\mu_\Om = \mu,
\end{equation}
where $\mu_\Om = (\g_\Om)_*\eta$ for some prescribed measure $\eta$ on $\del\Om$ and $\g_\Om$ is the Gauss map. 
The prescribed measure $\eta$ is typically absolutely continous with respect to $\h^{n-1}\on_{\del\Om}$ and with respect to the Gauss map, we can 
express the density of the induced measure on $\Snn$ as 
\begin{equation}\label{eq:msrdensity}
d\mu_\Om( \xi) = \F[h_\Om]( \xi) d \xi 
\end{equation}
for a functional $\F: \K(\Snn) \to \R$. Existence of domains for the Minkowski problem involves computation of the first variation of measures 
corresponding to $\Om^t= \Om + t\Om_0$ with the support functions $h_{\Om^t}= h_\Om +t h_{\Om_0}$, 
\begin{equation}\label{eq:L}
\LL_\Om [v] := \frac{d}{dt}\Big|_{t=0} \F[h_{\Om^t}] = \frac{d}{dt}\Big|_{t=0} \F[h_\Om +t h_{\Om_0}]. 
\end{equation}
This is used to show the existence of domains using continuity method or variational techniques of constrained 
minimization problems. 

In case of the original Minkowski problem, we have $\eta= \h^{n-1}$ and 
$\mu_\Om= (\g_\Om)_*\h^{n-1}$. The case of the given measure being discrete was considered in \cite{M1, M2}, where the 
corresponding convex domains are polytopes. The continuous case has been shown in \cite{A1, A2,FJ} and the smooth case in \cite{CY}. 
If $\Om$ is of class $C^2_+$ and the covariant derivatives are as in \eqref{eq:cov} defined by the charts as shown in the 
previous subsection,
then recalling \eqref{eq:detkinv} and \eqref{eq:intbd}, we note that 
the density \eqref{eq:msrdensity} in this case, is the reciprocal of the Gauss curvature, i.e. 
$$ d\mu_\Om( \xi)=\det(\covtwo  h_\Om + h_\Om \I) d \xi= \frac{d \xi}{\kr(\inv{\g_\Om}( \xi))}.$$
Furthermore, if $\mu_\Om = (\g_\Om)_*\eta$ where $d\eta = f d\h^{n-1}$ for a function $f:\del\Om \to \R$, then we have 
\begin{equation}\label{eq:fhom}
\F[h_\Om]( \xi) = \frac{(f\circ \inv{\g_\Om})( \xi)}{\kr(\inv{\g_\Om}( \xi))} = f\big(\inv{\g_\Om}( \xi)\big)\det(\covtwo  h_\Om + h_\Om \I).
\end{equation}
The above formulae hold for $S_\Om$-a.e. $ \xi\in \Snn$ for a general convex domain. 
As examples of the above, the prescribed measure is the harmonic measure $\om$ at origin in \cite{J1} where we have $\mu_\Om = (\g_\Om)_*\om$ and 
$f = (\del G/\del \nu)$ where $G$ is the Green's function with pole at $0$; in the case of capacitary measures in \cite{CNSXYZ}, we have 
$f= |\gr U|^p$ where $U$ is the capacitary function. 
In our case, we consider $\mu_\Om =  (\g_\Om)_*\om_p$ where $\om_p$ a the $p$ harmonic measure with respect to a function $u$ and the domain $\Om$ being convex, we have $f= |\gr u|^{p-1}$ for $1<p<\infty$.

\section{Existence for $p$-Harmonic Minkowski problem}\label{sec:ex}
In this section, we consider the Minkowski problem \eqref{eq:minkprob} for $p$-harmonic measures. We establish the existence of solutions 
of the problem. 
Henceforth, we shall denote 
\begin{equation}\label{eq:pharmink}
\mu_\Om =  (\g_\Om)_*\om_p,
\end{equation}
where $\om_p$ is the $p$ harmonic measure with respect to a function $u=u_\Om\in W^{1,p}(\Om\cap N)$ given by 
$d\om_p=|\gr u|^{p-1} d\h^{n-1}\on_{\del\Om}$ where $u$ is $p$-harmonic in $\Om\cap N$ and satisfies 
\begin{equation}\label{eq:omdir}
 \begin{cases}
  \dv (|\gr u|^{p-2}\gr u)=0\ \ &\text{in}\ \Om\cap N;\\
   u > 0 \ \ &\text{in}\ \Om;\\
  u= 0 \ \ &\text{on}\ \del\Om,
 \end{cases}
\end{equation}
where $N$ is a neighbourhood of $\del\Om$; thus, $u\in W^{1,p}(N)$ upon zero extension. 
Up to possible reduction, the choice of $N$ is made so that 
$\gr u\neq 0$ in $\Om\cap N$ and 
\begin{equation}\label{eq:udeln}
\|u\|_{L^\infty(\bar N\cap \Om)}+\|\gr u\|_{L^\infty(\bar N\cap \Om)} <\infty,
\end{equation}
and we also assume that $\del N$ is $C^\infty$. 
For a general bounded convex set $K$, if $K=\bar\Om$ for a convex domain $\Om$ then we define $\mu_K=\mu_\Om$; if $K$ is of empty interior, i.e. $K=\del K$, 
we define $\mu_K= (\g_K)_* \om_p$ in the sense of \eqref{eq:gmap} and \eqref{eq:ginv}, where $\om_p$ is the $p$-harmonic measure with respect to a function $u=u_K\in W^{1,p}(N)$ with $N$ being a neighbourhood of $K$, given by 
$d\om_p=|\gr u|^{p-1} d\h^{n-1}\on_{K}$ where $u$ is $p$-harmonic in $ N$ and satisfies 
\begin{equation}\label{eq:omdiremp}
 \begin{cases}
  \dv (|\gr u|^{p-2}\gr u)=0\ \ &\text{in}\  N\mns K;\\
   u \geq 0 \ \ &\text{in}\ N;\\
  u= 0 \ \ &\text{on}\  K.
 \end{cases}
\end{equation}
Thus, in general, 
we have $\mu_K$ defined on $\Snn$ associated to $u=u_K\in W^{1,p}(N)$ as 
\begin{equation}\label{eq:mukint}
\mu_K(E) = \int_{\inv{\g_K}(E)} |\gr u|^{p-1} d\h^{n-1}, \quad \text{for any measurable}\ E\sub \Snn,
\end{equation}
where $\g_K$ and $\inv{\g_K}(E)$ are in the sense of \eqref{eq:gmap} and \eqref{eq:ginv}. If $K=\bar \Om$ for a convex domain $\Om$ then \eqref{eq:mukint} coincides with the standard push forward measure \eqref{eq:pharmink}. 

\subsection{First variation of density of $p$-harmonic measures}\label{subsec:fstvar}
Here, we consider $\Om$ as strongly convex domain of class $C^{2,\alpha}_+$ so that $\g_\Om: \del \Om  \to  \mathbb{S}^{n-1}$ is a diffeomorphism. 
Let its support function be $h_\Om$. Using the same notations as Section \ref{sec:prelim}, 
from \eqref{eq:pharmink} and 
\eqref{eq:phm}, we have $$d\mu_\Om = |\gr u(F_\Om( \xi))|^{p-1}d\h^{n-1}\on_{\del\Om}=  \F[h_\Om]( \xi) d  \xi, $$
where $ F_\Om( \xi) := \inv{\g_\Om}( \xi)= \gr h_\Om( \xi) $ and recalling \eqref{eq:fhom}, we have 
$$  \F[h_\Om]( \xi) = |\gr u(F_\Om( \xi))|^{p-1} \det(\covtwo  h_\Om + h_\Om \I). $$
Henceforth, we shall denote $h=h_\Om$ and $F=F_\Om$. Thus, from \eqref{eq:supg} and \eqref{eq:hgradh}, we have 
$ h( \xi)= h_\Om ( \xi)= \inp{ \xi}{\inv{\g_\Om}( \xi)} = \inp{x}{\g_\Om(x)}$ and 
$$F( \xi)= F_\Om ( \xi)= \inv{\g_\Om}( \xi)= \gr h ( \xi) .$$ 
As in Section \ref{sec:prelim}, we shall consider the orthonormal frame field $\{e^1,\ldots, e^{n-1}\}$ of $\Snn$ so that for any $ \xi\in \Snn$, the tangent space $T_ \xi(\Snn)$ is spanned by $\{e^i( \xi)\}$; furthermore, we denote covariant derivatives with respect to the frame as in \eqref{eq:cov} with respect to the coordinate charts as in Section \ref{sec:prelim}. 
We have 
$ \det(\covtwo  h_\Om + h_\Om \I) = 1/(\kr \circ \inv{\g_\Om})$ from \eqref{eq:detkinv} and \eqref{eq:intbd},  and recalling \eqref{eq:phm}, for any integrable function $f:\del\Om\to \R$, we have 
\begin{equation}\label{eq:intbdom}
\int_{\del \Om} f(x) \, d\om_p(x)= \int_{\Snn} f(\inv{\g_\Om}( \xi)) |\gr u(\inv{\g_\Om}( \xi))|^{p-1}\det\big(\covtwo  h_\Om( \xi) + h_\Om( \xi) \I\big) \, d \xi. 
\end{equation}

We consider a convex domain $\Om_0\subset \R^n$, also of class $C^{2,\alpha}_+$ to find the first variation for $p$-harmonic measures on the following family of domains, 
\begin{equation}\label{eq:omt}
\Om^t= \Om+t\Om_0 \quad \text{with support functions}\quad h_{\Om^t}= h_\Om +t h_{\Om_0}=h+tv.
\end{equation}
Without loss of generality, we assume $0\in \Om_0$ so that $v\geq 0$ and hence we have $\Om \sub \Om^t$ for $t\geq 0$ and $\Om^t\sub \Om$ for $t\leq  0$, more generally 
$\Om^s\sub\Om^t$ if $s\leq t$. We take a small enough 
\begin{equation}\label{eq:tau}
\tau =\tau\big(n, d_\h(\del\Om,\del N), d_\h(\del\Om_0,\del N), \|u\|_{W^{1,p}(N)}\big)>0, 
\end{equation}
so that for all $|t|\leq\tau$ the domain $\Om^t$ is also of class $C^{2,\alpha}_+$ and $\del\Om^t \subset N$ and 
$\del N\cap \Om =\del N\cap \Om^{t}$. Also, note that, from \eqref{eq:hdh} we have $d_\h (\Om^t,\Om)= |t|\|v\|_{L^\infty(\Snn)}$ and hence 
$d_\h (\Om^t,\Om)\to 0^+$ as $t\to 0$. Moreover, from Steiner's formula for volumes and mixed volumes, we have
\begin{equation}\label{eq:st}
t\mapsto |\Om+t \Om_0|,\ \h^{n-1}(\del \Om+t \del \Om_0) \quad \text{are smooth;}
\end{equation}
in fact, the above maps are polynomials, see \cite[Chapter 5]{Sc}. 

 We consider the $p$-harmonic measures corresponding to $u(\cdot, t)\in W^{1,p}(\Om^t\cap N)$ defined as the weak solution of the following Dirichlet problem
\begin{equation}\label{eq:omtdir}
 \begin{cases}
  \dv \big(|\gr u(\cdot, t)|^{p-2}\gr u(\cdot, t)\big)=0,\ \ &\text{in}\ \Om^t\cap N;\\
  u(x, t)= 0, \ \  &\forall\ x\in \del\Om^t\cap N;\\
  u(x, t)= u\big(\frac{x}{1+t}\big),\ \  &\forall\ x\in \del N\cap \Om^t;
 \end{cases}
\end{equation}
for $t\in [-\tau, \tau]$ for $\tau>0$ small enough, so that upon zero extension, $u(\cdot, t)\in W^{1,p}(N)$. Then the corresponding measure is defined by 
$$d\mu_{\Om^t} = |\gr u (F( \xi, t), t)|^{p-1}d\h^{n-1}\on_{\del\Om^t}=  \F[h_{\Om^t}]( \xi) d  \xi$$ where 
$F( \xi, t):= F_{\Om^t}( \xi) = \inv{\g_{\Om^t}}( \xi)= \gr h( \xi) + t \gr v( \xi) $ and 
\begin{equation}\label{eq:fht}
\F[h+tv]( \xi) = |\gr u (F( \xi, t), t )|^{p-1} \det\big(\covtwo  h + h \I + t(\covtwo  v + v\I)\big).
\end{equation}
We need to establish some properties of the first variation as in \eqref{eq:L}. The behavior with respect to dilations is provided in the following lemma.
\begin{lemma}\label{lem:flhh}
Let $\F$ be as in \eqref{eq:fht} and $\LL_h[v]= \frac{d}{dt}\big|_{t=0} \F[h +t v]$, then we have 
\begin{equation}\label{eq:flhh}
\F[(1+t)h] = (1+t)^{n-p}\F[h] \quad \text{and}\quad \LL_h[h] = (n-p) \F[h],
\end{equation}
for all $|t|\leq\tau$ with a small enough $\tau>0$ as in \eqref{eq:tau}. 
\end{lemma}
\begin{proof}
Let $\lambda = 1+t$ and consider the case 
of $\Om_0=\Om$ and $v=h$ of the above. Then, we have $\Om^t= \lambda \Om$ and $ F_{\Om^t}=\inv{\g_{\Om^t}}= \lambda \gr h$.  
Let us denote $u_\lambda (\cdot) = u(\cdot, \lambda-1)$ as the weak solution of the problem \eqref{eq:omtdir} for this case, i.e.  
\begin{equation}\label{eq:lamomdir}
 \begin{cases}
  \dv \big(|\gr u_\lambda |^{p-2}\gr u_\lambda \big)=0,\ \ &\text{in}\ \lambda\Om\cap N;\\
  u_\lambda(x)= 0, \ \  &\forall\ x\in \del\lambda\Om\cap N;\\
  u_\lambda (x)= u(x/\lambda),\ \  &\forall\ x\in \del N\cap \lambda\Om;
 \end{cases}
\end{equation}
for $|\lambda-1|\leq\tau$ being small enough.  
Using these in \eqref{eq:fht}, we can conclude 
\begin{equation*}
\F[\lambda h]( \xi) = |\gr u_\lambda (\lambda \gr h( \xi))|^{p-1} \lambda^{n-1}\det\big(\covtwo  h + h \I \big) 
= \bigg(\frac{|\gr u_\lambda (\lambda \gr h( \xi))|}{|\gr u (\gr h( \xi))|}\bigg)^{p-1}\lambda^{n-1}\F[h]( \xi).
\end{equation*}
Now, note that \eqref{eq:omdir} trivially implies that $ u(x/\lambda)$ is also a solution of \eqref{eq:lamomdir} and hence, from uniqueness 
$ u_\lambda (x)= u(x/\lambda)$ in $\lambda\Om$. Therefore, using $\gr u_\lambda (x)= \frac{1}{\lambda}\gr u(x/\lambda)$ on the above, 
we get \eqref{eq:flhh} (the second part of \eqref{eq:flhh} follows immediately from the first) and the proof is finished. 
\end{proof}
Towards the computation of the first variation, we need to establish some basic properties. 
From \eqref{eq:c1reg}, \eqref{eq:c1regbd}, we know that weak solutions $u(\cdot, t)\in W^{1,p}(\Om^t\cap N)$ of \eqref{eq:omtdir} are $C^{1,\beta}$ in $\Om^t$ 
for any $t\in [-\tau, \tau]$ and some $\beta=\beta(n,p)\in (0,1)$. 

First, we show uniform continuity of $t\mapsto u(\cdot, t)$. 
\begin{lemma}\label{lem:unifc}
For any $t\in [-\tau, \tau]$, if $u(\cdot, t)\in C^{1,\beta}_\loc(\Om^t\cap N)$ is the solution of the Dirichlet problem \eqref{eq:omtdir}, then $t\mapsto u(\cdot, t), \gr u(\cdot, t) $ are 
uniformly continuous on compact subsets of $N$.  
\end{lemma}
\begin{proof}
From \eqref{eq:omtdir}, note that for any $\eta\in W^{1,p}(N)$, we have 
\begin{equation}\label{eq:weakomt}
\begin{aligned}
\int_{\Om^t\cap  N}  &|\gr u(x, t)|^{p-2}\inp{\gr u(x, t)}{\gr \eta(x)} \dx \\
&=\int_{\del(\Om^t\cap N)} \eta(x) |\gr u(x, t)|^{p-2}\inp{\gr u(x, t)}{\nu_{\del(\Om^t\cap N)(x)}} \,d \h^{n-1}(x) \\
&=\frac{1}{|1+t|^{p-2}(1+t)}\int_{\Om^t\cap \del N} 
\eta(x) |\gr u({\scriptstyle \frac{x}{1+t}})|^{p-2}\inp{\gr u({\scriptstyle \frac{x}{1+t}})}{\nu_{\Om^t\cap \del N}(x)} \,d \h^{n-1}(x)\\
&\qquad\qquad +\int_{\del\Om^t} \eta(x) |\gr u(x, t)|^{p-2}\inp{\gr u(x, t)}{\g_{\Om^t}(x)} \,d \h^{n-1}(x)
\end{aligned}
\end{equation}
for any $t\in [-\tau, \tau]$. From boundary conditions of \eqref{eq:omtdir} and comparison principle, we can conclude $u(x,t)> 0$ in $\Om^t\cap N$, hence 
$\g_{\Om^t}(x)= -\gr u(x,t)/|\gr u(x,t)|$ for all $x\in \del\Om^t, \ |t|<\tau$ for a small enough $\tau$ in \eqref{eq:tau}. Using \eqref{eq:weakomt} for $t,s\in [-\tau, \tau]$ and taking their difference,
we obtain 
\begin{equation*}
\begin{aligned}
\int_{\Om^t\cap  N}  &|\gr u(x, t)|^{p-2}\inp{\gr u(x, t)}{\gr \eta(x)} \dx -\int_{\Om^s\cap  N}  |\gr u(x, s)|^{p-2}\inp{\gr u(x, s)}{\gr \eta(x)} \dx \\
&=\frac{1}{|1+t|^{p-2}(1+t)}\int_{\Om^t\cap \del N} 
\eta(x) |\gr u({\scriptstyle \frac{x}{1+t}})|^{p-2}\inp{\gr u({\scriptstyle \frac{x}{1+t}})}{\nu_{\Om^t\cap \del N}(x)} \,d \h^{n-1}(x)\\
&\qquad - \frac{1}{|1+s|^{p-2}(1+s)}\int_{\Om^s\cap \del N} 
\eta(x) |\gr u({\scriptstyle \frac{x}{1+s}})|^{p-2}\inp{\gr u({\scriptstyle \frac{x}{1+s}})}{\nu_{\Om^s\cap \del N}(x)} \,d \h^{n-1}(x)\\
&\qquad - \Big(\int_{\del\Om^t} \eta(x) |\gr u(x, t)|^{p-1} \,d \h^{n-1}(x)-\int_{\del\Om^s} \eta(x) |\gr u(x, s)|^{p-1} \,d \h^{n-1}(x)\Big)
\end{aligned}
\end{equation*}
Then taking $\eta (x)=u(x,t)-u(x,s)$ on the above 
and using \eqref{eq:pmon}, it is not difficult to see that we can obtain 
$$ \int_{\Om^{\min\{t,s\}}\cap  N}  |\gr u(x, t)-\gr u(x,s)|^p \dx 
\leq c M^p \om (|t-s|), $$ 
where $c=c(n,p, \tau)>0, M=\| u\|_{L^\infty(\Om\cap \del N)}+ \|\gr u\|_{L^\infty(\Om\cap \del N)}
+\max_{|t|\leq\tau}\|u(\cdot,t)\|_{W^{1,p}(N)}<\infty$ and $\om (|t-s|) \to 0^+$ as $t\to s$, since $\gr u$ is continuous in $\Om\cap N$ and $\gr u(\cdot, t)$ is continuous in $\Om^t\cap N$ from \eqref{eq:c1reg} for all $t$, and furthermore, 
$|\Om^t| -|\Om^s|\to 0$ and $\h^{n-1}(\del\Om^t) -\h^{n-1}(\del\Om^s)\to 0$ as $t\to s$ from \eqref{eq:st}. The continuity of $x\mapsto \gr u(x,t)$ for every $t$ also implies that all points are Lebesgue points and hence the above is 
enough to conclude the following pointwise convergences $$|u(x,t)-u(x,s)|+|\gr u(x, t)-\gr u(x,s)|\to 0^+\quad \text{as}\ t\to s.$$ From \eqref{eq:c1reg}, 
we also have 
the following local estimates for any 
$B\subset\subset N$ with $M$ as above and $c=c(n,p,\tau,\diam(B))>0$, given by 
\begin{align*}
| u(x, t)- u(y,t)| \leq c M |x-y|  , \quad
|\gr u(x, t)-\gr u(y,t)| \leq c M |x-y|^\beta  ,\quad\forall\ x,y\in B,
\end{align*}
for every $t\in [-\tau, \tau]$, that affirms equicontinuity of $\{u(\cdot, t), \gr u(\cdot, t)\}_{0\leq |t|\leq\tau}$. Uniform continuity follows from the Arzel\`a-Ascoli theorem and the proof is finished. 
\end{proof}
\begin{remark}\label{rem:ux0}
Notice that, $u(\cdot, t)$ being the solution of \eqref{eq:omtdir}, from uniqueness and \eqref{eq:omdir} we have $u(x,0)= u(x)$ and hence, from the Lemma \eqref{lem:unifc},  
$u(x,t)\to u(x)$ uniformly as $t\to 0$. Hence, $\gr u (\cdot, t)\neq 0$ in $\Om^t\cap N$ and 
$\|u(\cdot, t)\|_{L^\infty(\bar N\cap \Om^t)}+\|\gr u(\cdot, t)\|_{L^\infty(\bar N\cap \Om^t)} <\infty $ for all $|t|\leq\tau$. 
\end{remark}
Now we need to investigate the differentiability of $t\mapsto u(\cdot,t)$. This is shown in the following that involve some standard regularity results and adaptations of the techniques in Colesanti et al. \cite{CNSXYZ}. 
For any function $f(\cdot, t)$ related to \eqref{eq:omtdir}, let us denote 
$$ \dot f (x) =  \frac{\del}{\del t}\Big|_{t=0} f(x, t),$$
if the derivative exists. Thus, we have $\dot F( \xi)= \gr v( \xi)$ and in the following proposition, we find $\dot u$.
Note that if $t\mapsto u(\cdot, t)$ is differentiable, then  
$\del_t(\dv (|\gr u|^{p-2}\gr u))= \dv (A(\gr u) \gr  \del_t u)=0$, where we henceforth denote
\begin{equation}\label{eq:Au}
 A(\gr u) := |\gr u|^{p-2} \left(\I + (p-2)\frac{\gr u\otimes \gr u}{|\gr u|^2} \right),
\end{equation}
which is uniformly elliptic modulo the weight $|\gr u|^{p-2}$. Indeed, it is not hard to see that 
\begin{equation}\label{eq:aubds}
 \min\{1,p-1\}|\gr u|^{p-2}| \zeta|^2\leq \inp{A(\gr u)  \zeta}{ \zeta} \leq \max\{1,p-1\}|\gr u|^{p-2}| \zeta|^2, \quad \forall\  \zeta\in \R^n. 
\end{equation}
However, the differentiability of $t\mapsto u(\cdot, t)$ being unknown apriori, difference quotients are to be used. This is shown in the following.  

\begin{proposition}\label{prop:udot}
If $u(\cdot, t)\in W^{1,p}(\Om^t\cap N)$ is the solution of \eqref{eq:omtdir} and $A$ is as in \eqref{eq:Au}, then the following holds: 
\begin{enumerate}
\item $t\mapsto u(\cdot, t)$ is differentiable at $t=0$ for all $x\in \bar \Om\cap N$ and $\dot u\in C^{2,\beta}(\br{\Om \cap N})$;
\item $\dot u$ is a solution of the equation $ \dv (A(\gr u) \gr \dot u)=0$ in $\Om\cap N$;
\item $\dot u(x) =-\inp{\gr u(x)}{x}$ for all $x\in \del N\cap \Om$;
\item $\dot u(x) = |\gr u (x)| v(\g_\Om(x))$ for all $x\in \del\Om$. 
\end{enumerate}
\end{proposition} 
\begin{proof}
For $t\neq 0$, let $ w(x,t) = (u(x,t)-u(x,0))/t$ and for any $s\in [0,1]$, let us denote
$$u_s(x,t)= su(x,t)+(1-s) u(x,0).$$ Then, $u_s':= d u_s/ds = t w(x,t)$ and for any $|t|\leq\tau$, note that we have 
\begin{equation}
\begin{aligned}
0&= \lap_p u(x,t)-\lap_p u(x,0) = \int_0^1 \frac{d}{ds} (\lap_p u_s (x,t) ) \, ds\\
&=  \int_0^1  \dv \big(A(\gr u_s(x,t)) \gr u_s'(x,t)\big) \, ds = t  \int_0^1  \dv \big(A(\gr u_s(x,t)) \gr w(x,t)\big) \, ds
\end{aligned}
\end{equation}
where $A$ is as in \eqref{eq:Au}. Thus, for $t\neq 0$, we have that $w(x,t)$ is a solution of 
the equation $\dv (B_u(x, t) \gr w(x, t)=0$, where entries of $B_u$ are given by 
\begin{equation}\label{eq:Bu}
B_u(x,t)_{i,j}=\inp{B_u(x,t)e_j}{e_i} = \int_0^1  \inp{A(\gr u_s(x,t)) e_j}{e_i} \, ds. 
\end{equation}
From Lemma \ref{lem:unifc} and Remark \ref{rem:ux0}, note that $\gr u_s(x,t)\to \gr u(x,0)=\gr u(x)$ uniformly as $t\to 0$, which together 
with \eqref{eq:Bu} implies $B_u(\cdot, t) \to A(\gr u)$ uniformly as well. Also, from Remark \ref{rem:ux0}, $\gr u(\cdot, t) \neq 0$ in $\Om^t\cap N$ for all $|t|\leq\tau$.  
Therefore, there exists $0<\lambda\leq\Lambda<\infty$ possibly dependent on 
$n,p, \min_{|t|\leq \tau} \inf_{x\in \Om^t\cap N}|\gr u(x,t)| $ and $\max_{|t|\leq \tau} \sup_{x\in \Om^t\cap N}|\gr u(x,t)|$ 
such that from \eqref{eq:aubds} and \eqref{eq:Bu}, we can conclude 
$$  \lambda | y|^2\leq \inp{B_u(x,t) y}{ y} \leq \Lambda | y|^2, \quad \forall\ x\in (\Om\cap\Om^t)\cap N,\  y\in \R^n,$$
for any $ t\in [-\tau,\tau]$, thereby making the equations uniformly elliptic. Hence, for $t\neq 0$, we have $w(\cdot, t)\in C^\infty (\Om\cap\Om^t\cap N)$. 
The limit of $w(x,t)$ needs to be established as $t\to 0$ 
for the proof. 

To this end, note that for $B\subset\subset \Om\cap N$ there exists $\tau'=\tau'(B)\leq \tau$ such that $B\cap \Om^t\neq \emp$ for 
$ t\in [-\tau',\tau']$. From \eqref{eq:c1reg} used on the equation \eqref{eq:omtdir} implies that $u(\cdot, t)\in C^{1,\beta}_\loc(\Om^t\cap N)$ for $\beta=\beta(n,p)\in (0,1)$, hence \eqref{eq:aubds} and \eqref{eq:Bu} gives $\|B_u(\cdot, t)_{i,j}\|_{C^{0,\sigma}(B\cap \Om^t)} \leq C'(n,p, M, \diam(B))$ 
for some $\sigma=\sigma(n,p,\beta)\in (0,1)$, where $M =\max_{|t|\leq\tau}\|u(\cdot,t)\|_{W^{1,p}(N)}$. Then, from a standard compactness argument 
$$ \|B_u(\cdot, t)_{i,j}\|_{C^{0,\sigma}(B)} \leq C(n,p, M,\diam(B)), \quad \forall\  t\in [-\tau',\tau'].$$
Therefore, we can invoke interior Schauder estimates, see \cite[Theorem 6.2]{GT}, for the equation 
$\dv (B_u(x, t) \gr w(x, t))=0$ on $ \Om\cap N$ to conclude that there exists $0<\tau'=\tau'(B)\leq \tau$ for any $B\subset\subset \Om\cap N$ such that 
\begin{equation}\label{eq:wc2a}
\|w(\cdot, t)\|_{C^{2,\beta}(B)} \leq \Lambda', \qquad \forall\ 0<|t|\leq \tau',
\end{equation}
for some $\Lambda' = \Lambda'\big(n,p, \lambda, M, \diam(B)\big)>0$ and $\beta\in (0,1)$. It is not hard to see that 
\eqref{eq:wc2a} implies uniform boundedness and equicontinuity of the family $\{w(\cdot, t)\}_{0<|t|\leq\tau'}$ and hence Arzel\`a-Ascoli theorem yields the existence of a convergent subsequence. For a sequence $t_k\to 0$ as $k\to \infty$, let $w\in W^{1,p}(\Om \cap N)$ be such 
a subsequential limit, i.e. 
\begin{equation}\label{eq:w}
w(x):= \lim_{k\to \infty} w(x, t_k)=  \lim_{k\to \infty} \frac{u(x,t_k)-u(x,0)}{t_k},\qquad \forall\ x\in \Om\cap N;
\end{equation}
then taking $t_k\to 0$ on \eqref{eq:wc2a} we get $w\in C^{2,\alpha}_\loc(\Om\cap N)$ and since $B_u(\cdot, t_k)\to A(\gr u)$ uniformly, we have that 
$w$ is a solution of $\dv (A(\gr u)\gr w)=0$ in $\Om\cap N$. 

Now we look into the boundary behavior of $w$. Since $\Om^t$'s are of class $C^{2,\alpha}$ and $u\in C^\infty (\Om^t\cap N)$ for all $|t|\leq \tau$, 
\eqref{eq:c1regbd} used for the equation \eqref{eq:omtdir} implies that $u(\cdot, t)\in C^{1,\gamma}(\br{\Om^t\cap N})$ for some $\gamma =\gamma(n,p,\alpha)\in (0,1)$
and $\|u(\cdot, t)\|_{C^{1,\gamma}(\Om^t\cap N)} \leq C(n,p, \alpha, M, N)$ with $M>0$ as above. This, together with \eqref{eq:aubds}, \eqref{eq:Bu} and 
\eqref{eq:st}, further imply
$$ \|B_u(\cdot, t)_{i,j}\|_{C^{0,\sigma}(\Om\cap N)} \leq C(n,p,\alpha,M, N), \quad \forall\  t\in [-\tau,\tau],$$
for some $\sigma=\sigma (n,p,\gamma)\in (0,1)$. Recalling
$\del N\cap \Om^t= \del N \cap \Om$ for all $|t|\leq \tau$ and from boundary condition of \eqref{eq:omtdir}, note that we have 
\begin{equation}\label{eq:wxtdn}
w(x,t) =\frac{u( \frac{x}{1+t})-u(x)}{t}= -\frac{\inp{\gr u(x)}{x}}{(1+t)}+ O(|t|^\gamma),\quad\forall\ x\in \del N\cap \Om^t= \del N \cap \Om, t\neq 0,
\end{equation}
which implies $\|w(\cdot, t)\|_{L^\infty(\del N\cap \Om)} \leq C(n,p,\tau, \gamma, M, N)$. Also from boundary condition of \eqref{eq:omtdir}, we can conclude that for some 
$\beta=\beta(n,p,\gamma)\in (0,1)$, 
\begin{equation}\label{eq:wc2ad}
\|w(\cdot, t)\|_{C^{2,\beta}(\del(\Om\cap \Om^t\cap N))} \leq C(n,p,\tau, \gamma, M, N), \quad\forall\ 0<|t|\leq \tau. 
\end{equation}
Therefore, we can invoke global Schauder estimates up to the boundary, see \cite[Theorem 6.6]{GT}, for the equation 
$\dv (B_u(x, t) \gr w(x,t))=0$ on $ (\Om\cap \Om^t)\cap N$ to conclude 
\begin{equation}\label{eq:wc2abd}
\|w(\cdot, t)\|_{C^{2,\beta}(\Om\cap \Om^t\cap N)} \leq \Lambda'', \qquad \forall\ 0<|t|\leq \tau,
\end{equation}
for some $\Lambda'' = \Lambda''\big(n,p,\alpha, \lambda,\tau, M, N\big)>0$ and $\beta\in (0,1)$. Taking $t_k\to 0$ on \eqref{eq:wc2ad} and \eqref{eq:wc2abd} and using \eqref{eq:st}, we get 
$w\in C^{2,\beta}(\br{\Om\cap N})$ and $\|w\|_{C^{2,\beta}(\Om\cap N)} \leq \Lambda''$. Using this we can find the boundary values of $w$ at $x\in \del(\Om\cap N)$ 
by taking sequences $x_j\in \Om\cap N$ such that $x_j\to x$ as $j\to \infty$ and using \eqref{eq:w} as
$$ w(x) = \lim_{j\to \infty} w(x_j) =\lim_{j\to \infty}  \lim_{k\to \infty} w(x_j, t_k)= \lim_{j\to \infty} \lim_{k\to \infty} \frac{u(x_j,t_k)-u(x_j,0)}{t_k}. $$
Since we have $\del N\cap \Om^{t_k}= \del N \cap \Om$, the evaluation of the limit is clear on $\del N\cap \Om$ taking $t_k\to 0$ and using
\eqref{eq:wxtdn} and continuity of $\gr u$, to obtain 
$$w(x)= -\inp{\gr u(x)}{x},\qquad \forall\ x\in \del N \cap \Om.$$ 
However, finding $w$ at $\del \Om$ is more involved since values at
$\del \Om^{t_k}$'s can accumulate to that on $\del\Om$ as $t_k\to 0$ and the limits may not be interchangeable in general. Here, the strong convexity of the domain is used to obtain sequences $x_k\in\del\Om^{t_k}$ with a common normal to evaluate the limit diagonally.  
Indeed, since $\Om^t$'s are of class $C^{2,\alpha}_+$, hence $\g_{\Om^t}:\del\Om^t\to \Snn$ is a diffeomorphism for every $|t|\leq \tau$. Hence, for any $x\in \del \Om$
with $\g_\Om (x)= \xi\in \Snn$, let $x_k= \inv{\g_{\Om^{t_k}}} ( \xi)\in \del\Om^{t_k}$. Then, recalling \eqref{eq:omt} and \eqref{eq:hgrh}, notice that 
\begin{equation}\label{xkx}
x_k =  \gr h( \xi) + t_k \gr v( \xi)= x + t_k \gr v ( \xi),
\end{equation}
and thus $x_k\to x$ as $k\to \infty$. Since \eqref{eq:wc2ad} and \eqref{eq:wc2abd} together with \eqref{xkx} imply 
$$|w(x_k,t_k)-w(x,t_k)|\leq C |x_k-x|=C|\gr v ( \xi)||t_k|,$$ 
for $C=C\big(n,p,\alpha, \lambda,\tau, M, N\big)>0$.
This is used in \eqref{eq:w} to obtain 
\begin{equation}\label{w1}
w(x) = \lim_{k\to \infty} w(x_k, t_k)= \lim_{k\to \infty} \frac{u(x_k,t_k)-u(x_k,0)}{t_k}=\lim_{k\to \infty}-\frac{u(x_k)}{t_k},
\end{equation}
since $u(x_k,t_k)=0$ from boundary conditions of \eqref{eq:omtdir} as $x_k\in \del\Om^{t_k}$. 
Furthermore, since $x\in \del \Om$, from boundary conditions of \eqref{eq:omdir}
we also have $u(x)=0$. Also, since $u\in C^{1,\gamma}(\br{\Om\cap N})$ with 
$\|u\|_{C^{1,\gamma}(\Om\cap N)} \leq C(n,p, \alpha, \|u\|_{L^\infty(\bar N\cap \Om)}+\|\gr u\|_{L^\infty(\bar N\cap \Om)}, N)$, therefore, 
using \eqref{xkx}, note that
$$ u(x_k)-u(x)= \inp{\gr u(x)}{x_k-x} + O(|x_k-x|^{1+\gamma})= t_k \inp{\gr u(x)}{\gr v ( \xi)} + O(|t_k|^{1+\gamma}).$$ 
Using these on \eqref{w1} we obtain
\begin{equation}\label{w2}
w(x) = \lim_{k\to \infty} \frac{u(x)-u(x_k)}{t_k} = -\inp{\gr u(x)}{\gr v( \xi)}, \qquad \forall\ x\in \del\Om. 
\end{equation}
Recall that, since $\del\Om= \del\{u>0\}$, the outer normal  $\nu_{\del\Om}= \xi= -\gr u/|\gr u| $ and 
since $v=h_{\Om_0}$, from \eqref{eq:hgrh}, $v( \xi)= \inp{ \xi}{\gr v( \xi)}$. Using these on \eqref{w2}, 
we finally obtain $w(x)=|\gr u (x)| v( \xi)$. 

Thus, we established that any subsequential limit $w$ as in \eqref{eq:w} must be a solution of the equation 
$\dv (A(\gr u)\gr w)=0$ in $\Om\cap N$ and moreover, $w(x)= -\inp{\gr u(x)}{x}$ for all $x\in \del N \cap \Om$ and 
$w(x)=|\gr u (x)| v(\g_\Om(x))$ for all $x\in \del\Om\cap N$. This equation being uniformly elliptic in $\Om\cap N$, 
$w$ is unique for every subsequence leading to $ \liminf_{t\to 0} w(\cdot, t)=w= \limsup_{t\to 0} w(\cdot, t)$. Therefore, the limit exists 
and hence $w=\dot u$, which completes the proof. 
\end{proof}
In view of the above lemma, let us define the operator $\Lambda_u: C(\del\Om)\to C(\del\Om)$ as 
\begin{equation}\label{eq:lambu}
\Lambda_u [f](x)= \inp{\gr w(x)}{\g_\Om(x)}, 
\end{equation}
for all $x\in \del\Om, f\in C(\del\Om)$, where 
$w\in  C^{1}(\br{\Om\cap N})$ is the solution of the Dirichlet problem 
\begin{equation}\label{eq:wdir}
 \begin{cases}
  \dv (A(\gr u) \gr w)=0\ \ &\text{in}\ \Om\cap N;\\
   w(x)=  -\inp{\gr u(x)}{x} \ \ &\forall\ x\in \del N\cap \Om;\\
   w(x)= f(x) \ \ &\forall\ x\in \del\Om\cap N;\\
 \end{cases}
\end{equation}
with $A(\gr u)$ as in \eqref{eq:Au}. 
Notice that Proposition \ref{prop:udot} implies that 
\begin{equation}\label{eq:udotlamb}
\Lambda_u \big[|\gr u| (v\circ \g_\Om)\big](x)=  \inp{\gr \dot u(x)}{\g_\Om(x)},\qquad\forall\ x\in \del\Om. 
\end{equation}
We have the following technical lemma. 
\begin{lemma}\label{lem:weta}
Let $w\in C^{1}(\br{\Om\cap N})$ be the solution of \eqref{eq:wdir}. For any $\eta\in C^1(\Om)$,  we have 
\begin{equation}\label{eq:weta}
\begin{aligned}
(p-1)\int_{\del\Om} \eta &|\gr u|^{p-2}  \inp{\gr w(x)}{\g_\Om(x)} d\h^{n-1}(x)\\
&= \int_{\Om\cap  N} |\gr u|^{p-2} \Big( \inp{\gr w}{\gr \eta} + \frac{(p-2)}{|\gr u|^2}  \inp{\gr w}{\gr u} \inp{\gr u}{\gr \eta} \Big)\dx\\
&\qquad\qquad + \int_{\del N\cap\Om} \eta  \,\inp{A(\gr u) \big(\gr u+D^2u\, x\big)}{\nu_{\del N\cap\Om}(x)} d\h^{n-1}(x). 
\end{aligned}
\end{equation}
\begin{proof}
If $w\in C^{1}(\br{\Om\cap N})$ is the solution of \eqref{eq:wdir}, then for any $\eta\in C^1(\Om)$, we have 
\begin{equation}\label{eq:wdirid}
\int_{\Om\cap N} \inp{A(\gr u) \gr w}{\gr \eta}\dx = \int_{\del(\Om\cap N)} \eta  \,\inp{A(\gr u) \gr w}{\nu_{\del (N\cap\Om)}(x)} d\h^{n-1}(x).
\end{equation}
Recalling \eqref{eq:Au}, we note that 
\begin{equation}\label{auw}
A(\gr u) \gr w = |\gr u|^{p-2} \Big(\gr w + \frac{(p-2)}{|\gr u|^2} \inp{\gr w}{\gr u} \gr u \Big).
\end{equation}
Since $\g_\Om(x)= -\gr u(x)/|\gr u(x)|$ for $x\in \del\Om$, hence the right hand side of \eqref{eq:wdirid}, together with \eqref{auw}, 
becomes 
\begin{align*}
\int_{\del\Om} \eta  &\,\inp{A(\gr u) \gr w}{\g_\Om(x)} d\h^{n-1}(x) + 
\int_{\del N\cap\Om} \eta  \,\inp{A(\gr u) \gr w}{\nu_{\del N\cap\Om}(x)} d\h^{n-1}(x)\\
&= \int_{\del\Om} \eta  |\gr u|^{p-2}\Big( \inp{\gr w}{\g_\Om(x)} + \frac{(p-2)}{|\gr u|^2} \inp{\gr w}{\gr u} \inp{\gr u}{\g_\Om(x)} \Big) d\h^{n-1}(x)\\
&\qquad \qquad-  \int_{\del N\cap\Om} \eta  \,\inp{A(\gr u) \big(\gr u+D^2u\, x\big)}{\nu_{\del N\cap\Om}(x)} d\h^{n-1}(x)\\
&=(p-1)\int_{\del\Om} \eta |\gr u|^{p-2}  \inp{\gr w(x)}{\g_\Om(x)} \,d\h^{n-1}(x)\\
&\qquad \qquad-  \int_{\del N\cap\Om} \eta  \,\inp{A(\gr u) \big(\gr u+D^2u\, x\big)}{\nu_{\del N\cap\Om}(x)} d\h^{n-1}(x). 
\end{align*}
Using \eqref{auw} on the left-hand side of \eqref{eq:wdirid} together with the above, we obtain \eqref{eq:weta} and the proof is finished. 
\end{proof}
\end{lemma}
Using the above, we show that $\Lambda_u$ is self-adjoint on $L^2(\del\Om, |\gr u|^{p-2}\h^{n-1})$ in the following.
\begin{proposition}\label{prop:selfadj}
Let $\Lambda_u: C(\del\Om)\to C(\del\Om)$ be as in \eqref{eq:lambu}. 
Then, for any $\psi_1, \psi_2 \in C(\del\Om)$, the following holds, 
\begin{equation}\label{eq:selfadj}
\int_{\del\Om} \psi_1 \Lambda_u [\psi_2] \,|\gr u|^{p-2} d\h^{n-1} = \int_{\del\Om} \psi_2 \Lambda_u [\psi_1] \,|\gr u|^{p-2} d\h^{n-1}. 
\end{equation}
\end{proposition}
\begin{proof}
Given any $\psi_1, \psi_2 \in C(\del\Om)$, there exists $w_1,w_2\in C^{1}(\br{\Om\cap N})$ that are solutions of \eqref{eq:wdir} with 
$w_1=\psi_1$ and $w_2=\psi_2$ on $\del\Om$ and $w_1=-\inp{\gr u}{x}=w_2$ on $\del N \cap\Om$. Hence, using \eqref{eq:weta} of Lemma \ref{lem:weta}, note that
\begin{align*}
\int_{\del\Om} \psi_1 \Lambda_u [\psi_2] &|\gr u|^{p-2} d\h^{n-1} \\
&=\int_{\del\Om} w_1  \inp{\gr w_2(x)}{\g_\Om(x)} |\gr u|^{p-2} d\h^{n-1}(x) \\
&= \frac{1}{p-1}\int_{\Om\cap N} |\gr u|^{p-2} \Big( \inp{\gr w_2}{\gr w_1} + \frac{(p-2)}{|\gr u|^2}  \inp{\gr w_2}{\gr u} \inp{\gr u}{\gr w_1} \Big)\dx\\
&\quad -  \frac{1}{p-1}\int_{\del N\cap\Om} \inp{\gr u}{x} \,\inp{A(\gr u) \big(\gr u+D^2u\, x\big)}{\nu_{\del N\cap\Om}(x)} d\h^{n-1}(x), 
\end{align*} 
which is symmetric with respect to $w_1$ and $w_2$. This completes the proof. 
\end{proof}
\begin{remark}
 For the special case $p=2$, notice that we have $A(\gr u)= \I$ and $\Lambda_u=\Lambda$ corresponds to the Neumann operator as in \cite{J1,J2}. 
\end{remark}

Our next goal is to compute the first variation explicitly to establish its self-adjointness. 
We shall denote 
$\cof[\covtwo  h + h \I]= \text{cofactor matrix of}\  (\covtwo  h + h \I)$ 
(transpose of adjoint) so that, 
\begin{equation}\label{eq:cofdet}
\cof[\covtwo  h + h \I] (\covtwo  h + h \I) = \det(\covtwo  h + h \I)\I = 1/(\kr \circ \inv{\g_\Om})\I.
\end{equation}
Furthermore, since $\inv{\mathcal W_K}= \covtwo  h + h \I$ for $K=\bar\Om$, note that the above can be rewritten as 
\begin{equation}\label{eq:cofwin}
\mathcal W_K(x) = \kr (x) \,\cof[\covtwo  h + h \I](\g_\Om(x)), \qquad \ \forall\ x\in \del\Om. 
\end{equation}

The function $u$ of \eqref{eq:omdir} being $p$-harmonic is $C^{1,\beta}$ up to the boundary from \eqref{eq:c1regbd}, is uniformly elliptic in a neighborhood of $\del \Om$ since $\Om$ is of class $C^{2,\alpha}_+$. Hence, from 
boundary Schauder estimates, see \cite{GT}, $D^2u$ is pointwise well-defined on $\del \Om$. 
The following lemma which prescribes the Hessian $D^2u$ on $\del\Om$, has been shown for the case $p=2$ in \cite[Lemma A]{J1}.
\begin{remark}\label{rem:sign}
The sign in the following Lemma 
is reversed to that of \cite[Lemma A]{J1} 
as the Green's function $G\leq 0$ in \cite{J1}, therefore our choice $u\geq 0$ corresponds to $u=-G$ for $p=2$.   
\end{remark}
\begin{lemma}\label{lem:grad2u} 
Let $u: \Om\to \R$ be as in \eqref{eq:omdir} and $\{e^1,\ldots, e^{n-1},\}$ be an orthonormal frame field of $\Snn$ such that 
for any $ \xi\in \Snn$ the unit vectors $e^i=e^i( \xi)\in \Snn$ span the tangent space $T_ \xi (\Snn)$. The covariant derivatives being defined as in \eqref{eq:cov} with respect to the suitable local coordinate charts related to the frame, 
we have the following:  
\begin{enumerate}
\item $\inp{D^2 u(F( \xi)) e^i}{e^j} = -\kr(F( \xi)) |\gr u(F( \xi))| \,\cof_{i,j}[\covtwo  h + h \I]$;
\item $\inp{D^2 u(F( \xi)) \xi}{e^i} = -\kr(F( \xi)) \sum_{j}\cof_{i,j}[\covtwo  h + h \I]\,\cov_j (|\gr u(F( \xi))|)$;
\item $\inp{D^2 u(F( \xi)) \xi}{ \xi}= \frac{1}{(p-1)}\kr(F( \xi)) |\gr u(F( \xi))|  \Tr (\cof[\covtwo  h + h \I])$;
\end{enumerate}
where $\cof_{i,j}[\cdot]= \inp{\cof[\cdot] e^j}{e^i}$ are entries of the cofactor matrix as in \eqref{eq:cofdet} 
for $i,j\in \{1,\ldots, n-1\}$ with respect to this frame 
and $F( \xi)= \inv{\g_\Om}( \xi)= \gr h( \xi) $. 
\end{lemma} 

\begin{proof}
Let $K=\bar\Om$ and the coordinate charts for the covariant derivatives be as in Section \ref{sec:prelim} in a neighborhood $U$ of 
$\xi \in \Snn$, i.e. for $0\in V\subset \R^{n-1}$ the coordinates $\chi :V \to U\subset  \Snn$ satisfying $\chi(0)=\xi$, $\del_j\chi (0)= e^j$ and 
$\del_{i,j}\chi(0)=-\delta_{i,j}\chi(0)$ for all $i,j\in \{1,\ldots, n-1\}$.

By differentiating $F\circ \chi= \inv{\g_\Om}\circ \chi= \gr h\circ \chi $, we note that $\del_j(F\circ \chi)= (D^2h \circ \chi)\del_j\chi$ 
for any $j\in \{1,\ldots, n-1\}$ and hence, at $\chi=\chi(z)\in U\subset \Snn$ for any $z\in V\subset \R^{n-1}$, we have
$\cov_j F( \chi)=D^2h( \chi)\del_j\chi$, in particular, at $0\in V$, we have  
\begin{equation}\label{f1}
\cov_j F( \xi)=D^2h( \xi)e^j= \inv{\mathcal W_K}(\xi)e^j= (\covtwo  h (\xi) + h (\xi)\I )e^j. 
\end{equation} 
Since any $ \chi\in U\subset \Snn$ is the outer unit normal at $F(\chi)\in \del\Om$ as $F= \inv{\g_\Om}$, hence we have 
$ \chi= -\gr u(F( \chi))/|\gr u(F( \chi))|$. Using this together with the above, 
we get
\begin{equation}\label{f2}
\inp{\cov_j F( \chi)}{\gr u(F( \chi))}= -|\gr u(F( \chi))| \inp{\cov_j F( \chi)}{ \chi}= -|\gr u(F( \chi))| \inp{D^2h( \chi)  \chi}{\del_j\chi}=0,
\end{equation}
since, recalling \eqref{eq:hgradhcov}, $D^2h( \chi)  \chi=0$. 
Now, we differentiate \eqref{f2} by covariant derivative $\cov_i$ at $\xi=\chi(0)$ for any $i\in \{1,\ldots, n-1\}$ to obtain 
\begin{equation}\label{f3}
\inp{D^2 u(F( \xi)) \cov_i F( \xi)}{\cov_j F( \xi)} + \inp{\cov_{i,j} F( \xi)}{\gr u(F( \xi))}=0.
\end{equation}
We can take $\chi$ defined on the lower hemisphere without loss of generality and 
$\cov_{i,j} F( \xi)$ can be computed transforming $z_i/\sqrt{\scriptstyle{1-\sum_{i=1}^{n-1}|z_i|^2}}\mapsto y_i$ so that 
$\chi(z)\mapsto (y_1,\ldots, y_{n-1},-1)/\sqrt{\scriptstyle{1-\sum_{i=1}^{n-1}|y_i|^2}}$, $ \chi(0)\mapsto (0,-1)$ and $e^i\mapsto e_i$ as in 
\cite{J1}. 
Then, using \eqref{f1} on \eqref{f3}, one can obtain 
\begin{equation}\label{f5}
\inp{D^2 u(F( \xi)) \inv{\mathcal W_K}( \xi) e^i}{\inv{\mathcal W_K}( \xi)  e^j} = -|\gr u(F( \xi))| \inp{\inv{\mathcal W_K}( \xi) e^i}{e^j},
\end{equation}
(see the end of \cite[p. 382]{J1}). 
We successively multiply the above entrywise with the entries of
$\mathcal W_K(F( \xi))=\kr(F( \xi))\, \cof[\covtwo  h + h \I]$ 
recalling \eqref{eq:cofwin}, then sum over to obtain $(1)$. 

Next, we differentiate the equation 
$|\gr u(F( \chi))|= -\inp{\gr u(F( \chi))}{ \chi}$ at $\xi=\chi(0)$ to obtain
$$ \cov_j(|\gr u(F( \xi))|) = -\inp{D^2 u(F( \xi)) \cov_j F( \xi)}{ \xi} -\inp{\gr u(F( \xi))}{e^j}=-\inp{D^2 u(F( \xi)) \inv{\mathcal W_K}( \xi) e^j}{ \xi},$$
as $\inp{\gr u(F( \xi))}{e^j}= -|\gr u(F( \xi))|\inp{ \xi}{e^j}=0$ from orthonormality. Multiplying with the entries of 
$\mathcal W_K(F( \xi))=\kr(F( \xi)) \,\cof[\covtwo  h + h \I]$ from \eqref{eq:cofwin} and summing over, we can similarly obtain $(2)$. 

Finally, since $u$ is a solution of the equation \eqref{eq:omdir}, we have 
$$0=\dv (|\gr u|^{p-2}\gr u)= |\gr u|^{p-2}\Big( \lap u+ \frac{(p-2)}{|\gr u|^2}\inp{D^2u \gr u}{\gr u}\Big), $$ 
in $\Om\cap N$ and $u\in C^\infty(\Om\cap N)$. Taking pointwise limit, together with $(1)$, we get
\begin{align*}
(p-2)\inp{D^2 u(F( \xi)) \xi}{ \xi}&= -\lap u(F( \xi)) = -\Tr (D^2 u(F( \xi)))\\
&= -\sum_i  \inp{D^2 u(F( \xi)) e^i}{e^i}-\inp{D^2 u(F( \xi)) \xi}{ \xi}\\
&= \kr(F( \xi)) |\gr u(F( \xi))|  \Tr (\cof[\covtwo  h + h \I])-\inp{D^2 u(F( \xi)) \xi}{ \xi},
\end{align*}
which yields $(3)$. Thus the proof is finished. 
\end{proof}

Now we obtain the first variation of the density of $p$-harmonic measures in the following which has been shown previously for a special case of $p=2$ by Jerison \cite[Proposition 2]{J1}. 
\begin{proposition}\label{prop:fstvar}
Given support functions $h=h_\Om$ and $v=h_{\Om_0}$, the first variation, defined as $\LL_h[v]= \frac{d}{dt}\big|_{t=0} \F[h +t v]$ as in \eqref{eq:L} with $\F$ as in \eqref{eq:fht}, is given by 
\begin{equation}\label{eq:llv}
\LL_h[v]= \sum_{i,j=1}^{n-1}\cov_j\Big(\cof_{i,j}[\covtwo  h + h \I]|\gr u|^{p-1} \cov_i v\Big)
-\frac{(p-1)}{\kr}|\gr u|^{p-2} \Lambda_u \big[|\gr u| (v\circ \g_\Om)\big],
\end{equation}
where $\cof_{i,j}[\cdot]= \inp{\cof[\cdot] e^j}{e^i}$ are entries of the cofactor matrix as in \eqref{eq:cofdet} 
for $i,j\in \{1,\ldots, n-1\}$ with respect to a local orthonormal frame $\{e^i\}$ and $\Lambda_u$ is as in \eqref{eq:lambu}.
\end{proposition}
\begin{proof}
It is standard as shown in \cite{CY}, that 
$$ \frac{d}{dt}\Big|_{t=0}\det\big(\covtwo  h + h \I + t(\covtwo  v + v\I)\big)= \Tr\big( \cof[\covtwo  h + h \I] (\covtwo  v + v \I)\big).$$
Recalling $\F[h+tv]( \xi) = |\gr u (F( \xi, t), t )|^{p-1} \det\big(\covtwo  h + h \I + t(\covtwo  v + v\I)\big)$ and 
using the above, note that
\begin{equation}\label{lv1}
\begin{aligned}
\LL_h[v]&= |\gr u(F( \xi))|^{p-1}\frac{d}{dt}\Big|_{t=0}\det\big(\covtwo  h + h \I + t(\covtwo  v + v\I)\big) \\
&\qquad \qquad+ \det(\covtwo  h + h \I)\frac{d}{dt}\Big|_{t=0}  |\gr u (F( \xi, t), t )|^{p-1} \\
&=  |\gr u(F( \xi))|^{p-1} \Tr\big( \cof[\covtwo  h + h \I] (\covtwo  v + v \I)\big)\\
&\qquad \qquad +(p-1)|\gr u(F( \xi))|^{p-2} \det(\covtwo  h + h \I)\frac{d}{dt}\Big|_{t=0}|\gr u (F( \xi, t), t )|,
\end{aligned}
\end{equation}
and we need to compute only the last term of the above. Notice that, from the boundary condition of \eqref{eq:omtdir},
the outer unit normal at $F( \xi, t)\in \del\Om^t$ is given by
$$ \xi=\g_{\Om^t}(F( \xi, t))= -\frac{\gr u (F( \xi, t), t)}{|\gr u(F( \xi, t), t)|} ,$$
hence, $|\gr u(F( \xi, t), t)|=-\inp{\gr u (F( \xi, t), t)}{ \xi}$. Therefore, we obtain 
\begin{equation}\label{lv2}
\begin{aligned}
\frac{d}{dt}\Big|_{t=0}|\gr u (F( \xi, t), t )|&= - \frac{d}{dt}\Big|_{t=0}\inp{\gr u (F( \xi, t), t)}{ \xi}\\
&= -\Big( \inp{D^2 u(F( \xi))\dot F( \xi)}{ \xi} + \inp{\gr \dot u(F( \xi))}{ \xi} \Big)= -(I_1+I_2).
\end{aligned}
\end{equation}
Since, $\dot F( \xi)= \gr v( \xi)$ and recalling \eqref{eq:grhcov}, we have 
$\gr v( \xi)=v( \xi) \xi+\sum_{i=1}^{n-1} \cov_i v( \xi) e^i$. Using this together with Lemma \ref{lem:grad2u}
we compute the first term of \eqref{lv2} as 
\begin{equation*}
\begin{aligned}
I_1&= \inp{D^2 u(F( \xi))  \xi}{ \xi} v( \xi) +\sum_{i=1}^{n-1} \inp{D^2 u(F( \xi))  \xi}{e^i} \cov_i v( \xi)\\
&= \frac{1}{(p-1)}\kr(F( \xi)) |\gr u(F( \xi))| v( \xi) \Tr (\cof[\covtwo  h + h \I]) \\
&\qquad\qquad-\kr(F( \xi)) \sum_{i,j=1}^{n-1}\cof_{i,j}[\covtwo  h + h \I]\,\cov_j (|\gr u(F( \xi))|) \cov_i v( \xi). 
\end{aligned}
\end{equation*}
It has been shown in \cite{CY} that $\sum_j \cov_j \,\cof_{i,j}[\covtwo  h + h \I]=0$ for any $i$ and using this on the latter term of the above, 
we obtain
\begin{equation}\label{lv2i1}
\begin{aligned}
I_1&=  -\kr(F( \xi))\sum_{i,j=1}^{n-1}\cov_j\big(\cof_{i,j}[\covtwo  h + h \I]|\gr u(F( \xi))| \big) \cov_i v( \xi)\\
&\qquad +\frac{1}{(p-1)}\kr(F( \xi))|\gr u(F( \xi))| v( \xi) \Tr (\cof[\covtwo  h + h \I]).
\end{aligned}
\end{equation}
Using \eqref{lv2i1} on \eqref{lv2}, we obtain 
\begin{equation}\label{lv2fin}
\begin{aligned}
\frac{d}{dt}\Big|_{t=0}|\gr u (F( \xi, t), t )| 
&= \kr(F( \xi)) \sum_{i,j=1}^{n-1}\cov_j\big(\cof_{i,j}[\covtwo  h + h \I]|\gr u(F( \xi))| \big) \cov_i v( \xi)\\
&\qquad -\frac{\kr(F( \xi))}{(p-1)}|\gr u(F( \xi))| v( \xi) \Tr (\cof[\covtwo  h + h \I])-\inp{\gr \dot u(F( \xi))}{ \xi},
\end{aligned}
\end{equation}
Recalling $\det(\covtwo  h + h \I)=1/\kr(F( \xi))$ from \eqref{eq:detkinv}, we use \eqref{lv2fin} in \eqref{lv1} to get
\begin{equation}\label{lv3fin}
\begin{aligned}
\LL_h[v]
&=  |\gr u(F( \xi))|^{p-1} \Tr\big( \cof[\covtwo  h + h \I] (\covtwo  v + v \I)\big)\\
&\qquad +(p-1)|\gr u(F( \xi))|^{p-2} \sum_{i,j=1}^{n-1}\cov_j\big(\cof_{i,j}[\covtwo  h + h \I]|\gr u(F( \xi))| \big) \cov_i v( \xi)\\
&\qquad -|\gr u(F( \xi))|^{p-1}  v( \xi) \Tr (\cof[\covtwo  h + h \I]) \\
&\qquad- (p-1)\frac{|\gr u(F( \xi))|^{p-2}}{\kr(F( \xi))} \inp{\gr \dot u(F( \xi))}{ \xi}.
\end{aligned}
\end{equation}
Now, we use the following identity to replace the second term of the above,
\begin{equation}\label{lvdiv}
\begin{aligned}
\sum_{i,j=1}^{n-1}\cov_j&\Big(\cof_{i,j}[\covtwo  h + h \I]|\gr u(F( \xi))|^{p-1} \cov_i v( \xi)\Big) \\
&= \sum_{i,j=1}^{n-1}\cov_j\Big(\cof_{i,j}[\covtwo  h + h \I]|\gr u(F( \xi))|^{p-1} \Big)\cov_i v( \xi) \\
&\quad\qquad+  \sum_{i,j=1}^{n-1} \cof_{i,j}[\covtwo  h + h \I]|\gr u(F( \xi))|^{p-1} \cov_{j,i} v( \xi) \\
&=(p-1)|\gr u(F( \xi))|^{p-2} \sum_{i,j=1}^{n-1}\cov_j\big(\cof_{i,j}[\covtwo  h + h \I]|\gr u(F( \xi))| \big) \cov_i v( \xi)\\
&\quad\qquad+  |\gr u(F( \xi))|^{p-1} \Tr\big( \cof[\covtwo  h + h \I] \covtwo  v\big).
\end{aligned}
\end{equation}
Using \eqref{lvdiv} on \eqref{lv3fin}, we observe that the first and third term in \eqref{lv3fin} get cancelled off and we are left with 
$$ \LL_h[v]= \sum_{i,j=1}^{n-1}\cov_j\Big(\cof_{i,j}[\covtwo  h + h \I]|\gr u(F( \xi))|^{p-1} \cov_i v( \xi)\Big)
-(p-1)\frac{|\gr u(F( \xi))|^{p-2}}{\kr(F( \xi))} \inp{\gr \dot u(F( \xi))}{ \xi},$$
which completes the proof since $\Lambda_u \big[|\gr u| (v\circ \g_\Om)\big](F( \xi))=  \inp{\gr \dot u(F( \xi))}{ \xi}$ from \eqref{eq:udotlamb}.  
\end{proof}
The following shows the self-adjointness of the first variation in $L^2(\Snn, \h^{n-1})$.  
\begin{corollary}\label{cor:llvsym}
The first variation $\LL_h$ being as above, for any $v_1, v_2\in C(\Snn)$, we have 
\begin{equation}\label{eq:llvsym}
\int_{\Snn} v_1 \,\LL_h [v_2]\, d \xi \,=\, \int_{\Snn} v_2\, \LL_h [v_1] \, d \xi.  
\end{equation}
\end{corollary}
\begin{proof}
From Proposition \ref{prop:fstvar} and \eqref{eq:detkinv}, we have
\begin{equation*}
\begin{aligned}
\int_{\Snn} v_1( \xi) \LL_h [v_2]( \xi)\, d \xi &= \int_{\Snn}\sum_{i,j=1}^{n-1}v_1( \xi)\cov_j
 \Big(\cof_{i,j}[\covtwo  h + h \I]|\gr u(F( \xi))|^{p-1} \cov_i v_2( \xi)\Big)\, d \xi \\
&\quad -(p-1)\int_{\Snn} v_1( \xi)\Lambda_u  \big[|\gr u| (v_2\circ \g_\Om)\big](F( \xi)) \frac{|\gr u(F( \xi))|^{p-2}}{\kr(F( \xi))} \, d \xi\\
&= J_1+ J_2. 
\end{aligned}
\end{equation*}
Note that, using integral by parts (Stoke's theorem for compact manifold without boundary) on the first term of the above, we have 
$$ J_1= -\int_{\Snn}\sum_{i,j=1}^{n-1}|\gr u(F( \xi))|^{p-1} \cof_{i,j}[\covtwo  h + h \I] \cov_i v_2( \xi)\cov_jv_1( \xi)\, d \xi $$
which is symmetric in $v_1$ and $v_2$. For the second term, using \eqref{eq:detkinv} and \eqref{eq:intbd}, we have 
$$ J_2= -(p-1)\int_{\del\Om} (v_1\circ \g_\Om) \Lambda_u \big[|\gr u| (v_2\circ \g_\Om)\big] |\gr u|^{p-2} \, d\h^{n-1}, $$
which is also symmetric in $v_1$ and $v_2$ from self-adjointness \eqref{eq:selfadj} of $\Lambda_u$ from Proposition \ref{prop:selfadj}. 
Combining both cases of the above, the proof is finished. 
\end{proof}

\subsection{Convergence of $p$-harmonic measures}
The goal here is to establish weak convergence of $p$-harmonic measures with respect to the Hausdorff distance. 
Given 
a convex domain $\Om\subset\R^n$, a neighborhood $N$ of $\del\Om$ and $u\in W^{1,p}(\Om\cap N)$ as in \eqref{eq:omdir}, 
we consider a sequence of convex domains $\Om_j\subset \R^n$ such that 
$d_\h(\Om_j,\Om)\to 0^+$ as $j\to \infty$; recalling \eqref{eq:hdh}, $d_\h(\del\Om_j,\del\Om)\to 0^+$. 

First, we have the following lemma due to Jerison \cite{J1}.
\begin{lemma}\label{lem:nomj}
Given $\Om\subset\R^n$ and 
$d_\h(\Om_j,\Om)\to 0^+$ as $j\to \infty$, 
if $0\in \Om\cap\Om_j$ for all $j\geq j_0$ for some $j_0\in \N$,  
then there exists $r>0$
such that, up to a subsequence, $(1-1/j)\Om_j/r \subset \Om/r$.
\end{lemma}
\begin{proof}
Note that $|\dist(0,\del\Om_j)-\dist(0,\del\Om)|\leq d_\h(\del\Om_j,\del\Om)=d_\h(\Om_j,\Om)\leq 1/j$, for all $j\geq j_0$ with $j_0\in \N$ large enough. Hence, we choose 
$$0<r< \min\{\dist(0,\del\Om), \inf_{j\geq j_0} \dist(0,\del\Om_j)\} $$  
so that $B_r(0)\subset \Om_j\cap \Om$ for all $j\geq j_0$. 
Let $\tilde \Om=\Om/r$ and $\tilde \Om_j= \Om_j/r$.
Hence $B_1(0)\subset \tilde \Om_j\cap \tilde \Om$ for all $j\geq j_0$ and $d_\h(\tilde\Om_j,\tilde\Om)\to 0^+$. 
Thus, $1=h_{B_1(0)}\leq \min\{h_{\tilde\Om_j}, h_{\tilde\Om}\}$ and 
recalling \eqref{eq:hdh}, we have the convergence of support functions. Hence, up to a subsequence, 
$ 1 \leq h_{\tilde\Om_j}\leq h_{\tilde\Om}+1/j$ 
for all $j\geq j_0$ which implies $(1-1/j) h_{\tilde\Om_j}\leq h_{\tilde\Om}- (h_{\tilde\Om}-1)/j -1/j^2\leq h_{\tilde\Om}$, 
and the proof is complete. 
\end{proof}

Note that $\lim_{j\to\infty}d_\h \big((1-1/j)\Om_j,\Om\big)=\lim_{j\to\infty}d_\h(\Om_j,\Om)=0$. 
Henceforth, in view of Lemma \ref{lem:nomj}, we shall assume without loss of generality that 
\begin{equation}\label{eq:omjcond}
0\in \Om\cap\Om_j\qquad \text{and}\qquad \Om_j\subset \Om, \qquad \forall\ j\geq j_0, 
\end{equation}
for a large enough $j_0\in\N$. We can also similarly assume without loss of generality that $\Om \subset \Om_j (1+1/j)$ for all $j\geq j_0$. 
Furthermore, since $d_\h(\del\Om_j,\del\Om)\to 0^+$, for any sequence of
neighborhoods $N_j$ of $\del\Om_j$ we have $d_\h(N_j,N)\to 0^+$ as $j\to \infty$. We can assume without loss 
of generality that $N$ is a convex ring (see \cite{Lewis77}, e.g. annulus) so that 
the outer and inner segments of $\del N$ enclose convex domains. Hence, in view of \eqref{eq:hdh}
and Lemma \ref{lem:nomj}, we can regard that being subject to a possible small extension and dilation, $N$ is also a neighborhood of 
$\del\Om_j$ so that $\del N\cap \Om_j=\del N\cap \Om$ for all $j\geq j_0$, which shall be assumed hereafter.  

There exists a unique $u_j\in W^{1,p}(\Om_j\cap N)$ which is the weak solution of 
the Dirichlet problem
\begin{equation}\label{eq:omjdir}
 \begin{cases}
  \dv \big(|\gr u_j|^{p-2}\gr u_j\big)=0,\ \ &\text{in}\ \Om_j\cap N;\\
  u_j=0, \ \  &\text{in}\ \del\Om_j\cap N;\\
  u_j= u,\ \  & \text{in}\ \del N\cap \Om_j;
 \end{cases}
\end{equation}
for all $j\geq j_0$ and upon zero extension, $u_j\in W^{1,p}(N)$. 
We consider the $p$-harmonic measure with respect to $u_j$ as $\om_{p,j}$ 
so that $d\om_{p,j}= |\gr u_j |^{p-1}d\h^{n-1}\on_{\del\Om_j}$.

\begin{remark}\label{rem:conn}
In view of the above, given $d_\h \big(\Om_j,\Om)\to 0^+$ as $j\to\infty$, we can assume without loss of generality 
that $\Om_j(1-1/j)\subset \Om\subset \Om_j(1+1/j)$ for $j\geq j_0$ large enough. Hence, \eqref{eq:omjdir} can be re-defined 
with dilations on $\Om$ as $\Om/(1+1/j)$ or $\Om/(1-1/j)$ for which points in $\Om$ are transformed as $x\mapsto x/(1+1/j)$ and $x\mapsto x/(1-1/j)$ respectively. 
In this sense, the approximants defined by \eqref{eq:omjdir} can be regarded as a generalization of \eqref{eq:omtdir} for $t>0$ or $t<0$, 
since for $\Om^t=\Om+t\Om_0$, recall that 
$d_\h (\Om^t,\Om)= |t|\|h_{\Om_0}\|_{L^\infty(\Snn)}\to 0^+$ as $t\to 0$. 
All the following results corresponding to \eqref{eq:omjdir} remain unaffected with respect to these dilations. 
\end{remark}

Towards the weak convergence of $p$-harmonic measures, we make use of the following geometric lemma due to Jerison \cite[Lemma 3.3]{J2}, which  says that most of the boundaries of convex domains close enough by Hausdorff distance, can be locally flattened, except for a small measure set. It can be shown using  
convexity and Vitali covering theorem.
\begin{lemma}\label{lem:jerison}
 For a convex domain $\Om$, given any $ \eps > 0 $ there exists $\delta=\delta(\eps)>0$ and 
 a finite collection of disjoint balls $ B_{r_k} ( x_k)$ for $1\leq k \leq N_0,$ with $ r_k < \eps$ and $x_k  \in \del \Om$, such that for every $\Om'$ with 
 $d_\h(\Om',\Om)<\delta$, we have 
\begin{align*}  
\mathcal{H}^{n-1}    \big(\del\Om  \mns \bigcup_{k=1}^{N_0}    \bar B_{r_k} (x_k) \big)  < \eps,
\end{align*} 
and after a possible rotation and translation, $\del\Om, \del\Om'$ are graphs of functions 
$\phi,\phi'\in C^{0,1}(\R^{n-1})$ in $B_{r_k/\eps}(x_k)$, satisfying $| \gr \phi (x)| +  | \gr \phi' (x)|  \leq \eps$ for all $x\in \R^n$ with $|x|<r_k/\eps$. 
\end{lemma} 
The following is a part of \cite[Lemma 3.7]{J2} and can be shown using the density of radial projection in \eqref{eq:intbdal}, i.e. 
$d\h^{n-1} = |\rho_\Om( \theta)|^n/h_\Om (\g_\Om(\rho_\Om( \theta))) \, d\theta$ on Lemma \ref{lem:jerison}. 
\begin{corollary}\label{cor:jerison}
For a convex domain $\Om$, given any $ \eps > 0 $ there exists $\delta=\delta(\eps)>0,\ s_0>0$ and a family of balls $\mathcal B$ on $\Snn$ such that the folllowing holds:
\begin{enumerate}
\item Every $B\in \mathcal B$ has radius $s_0$.
\item If $B_r(0)\subset \Om\subset B_R(0)$, then there exists $N_0>0$ depending on $R/r$ such that any point in $\Snn$ belongs to 
atmost $N_0$ balls of $\mathcal B$.
\item $\theta( \Snn\mns F)<\eps$ where $F= \bigcup_{B\in \mathcal B} B$ and $\theta=\h^{n-1}\on_{\Snn}$ is the uniform measure. 
\end{enumerate}
\end{corollary}
Now, we sketch the weak convergence of $p$-harmonic measures in the following. 
\begin{proposition}\label{prop:pconvmsr}
Given a bounded convex domain $\Om\subset\R^n$, for any sequence of convex domains $\Om_j\subset \R^n$ with
$d_\h(\Om_j,\Om)\to 0^+$ as $j\to \infty$, if $\om_{p,j}$ is the $p$-harmonic measure with respect to the solution 
$u_j\in W^{1,p}(N)$ of \eqref{eq:omjdir}, then up to a subsequence, 
$\gr u_j \to \gr u$ uniformly in $N$ and 
we have the weak convergence $\om_{p,j}\wto \om_p$ as $j\to \infty$, i.e. for any $\eta\in C(\bar N)$, we have 
\begin{equation}\label{eq:pconvmsr}
\lim_{j\to \infty}\int_{\del\Om_j} \eta |\gr u_j|^{p-1} d\h^{n-1} = \int_{\del\Om} \eta |\gr u|^{p-1} d\h^{n-1} . 
\end{equation}
\end{proposition}
\begin{proof}
From boundary conditions of \eqref{eq:omjdir} and comparison principle, we have $u_j>0$ in $\Om_j$ and hence, 
$\g_{\Om_j}= -\gr u_j/|\gr u_j|$. Then, testing the equation \eqref{eq:omjdir} with any $\eta\in W^{1,p}(N)$, we have 
\begin{equation*}
\begin{aligned}
\int_{\Om_j\cap  N}  |\gr u_j|^{p-2}\inp{\gr u_j}{\gr \eta} \dx 
&=\int_{\del(\Om_j\cap  N)} \eta |\gr u_j|^{p-2}\inp{\gr u_j}{\nu_{\del(\Om_j\cap  N)}} \,d \h^{n-1} \\
&=\int_{\Om_j\cap \del N} 
\eta |\gr u|^{p-2}\inp{\gr u}{\nu_{\Om_j\cap \del N}} \,d \h^{n-1}
+\int_{\del\Om_j} \eta |\gr u_j|^{p-1} \,d \h^{n-1}.
\end{aligned}
\end{equation*}
Since $\del N\cap \Om_j=\del N\cap \Om$ for all $j\geq j_0$, we test the equation \eqref{eq:omdir} with $\eta$ and take the difference with the above to obtain
\begin{equation}\label{eq:wcd}
\begin{aligned}
\int_{\Om_j\cap  N} |\gr u_j|^{p-2}&\inp{\gr u_j}{\gr \eta} \dx -\int_{\Om\cap  N}  |\gr u|^{p-2}\inp{\gr u}{\gr \eta} \dx \\
&= - \Big(\int_{\del\Om_j} \eta |\gr u_j|^{p-1} \,d \h^{n-1}-\int_{\del\Om} \eta |\gr u|^{p-1} \,d \h^{n-1}\Big).
\end{aligned}
\end{equation}

Now we show that $\gr u_j\to \gr u$ uniformly in $N$ as $j\to \infty$. However, in this case $|\Om\mns \Om_j| $ or $|\h^{n-1}(\del\Om)-\h^{n-1}(\del\Om_j)|$ are not generally 
dominated by $d_\h (\Om_j,\Om)$. Therefore, Lemma \ref{lem:jerison} and Corollary \ref{cor:jerison} can be used by taking radial projections 
of \eqref{eq:intbdal} and decompose the integral on $\Snn$ into $F= \bigcup_{B\in \mathcal B} B$ and $ \Snn\mns F$ as in Corollary \ref{cor:jerison}, then obtain 
for any $\eps>0$,  
 $$ \Big|\int_B |\gr (u_j\circ \rho_j)|^{p-1}-|\gr (u\circ \rho)|^{p-1} d\theta \Big| \leq c\eps \int_B |\gr (u\circ \rho)|^{p-1} d\theta,\quad \forall\ j\geq j_0(\eps),\ B\in \mathcal B, $$
where $\rho_j\to \rho$ $\h^{n-1}$-a.e. in $\Snn$. This involves refined estimates of $p$-harmonic functions in \cite{CNSXYZ, LN2, LN3}, etc. used to obtain e.g. \cite[Lemma 4.5]{CNSXYZ} to which we refer to and omit the details here. The above together with $\theta( \Snn\mns F)<\eps$ from Corollary \ref{cor:jerison} allows us to estimate the integral on whole of $\Snn$ and then using \eqref{eq:dradconv} and \eqref{eq:intbdal}, we have 
$$\Big|\int_{\del \Om_j} |\gr u_j|^{p-1}d\h^{n-1} -\int_{\del\Om}|\gr u|^{p-1} d\h^{n-1} \Big| \leq c\eps \|u\|_{W^{1,p}(N)}^{p-1},\quad \forall\ j\geq j_0(\eps).$$
An integral on $\Om\mns\Om_j$ can be similarly estimated using \eqref{eq:polark}, Corollary \ref{cor:jerison} since $r_{\Om_j}\to r_\Om$ uniformly as $j\to \infty$.  
Using these estimates together with $\eta= u_j-u$ on \eqref{eq:wcd} and using the monotonicity \eqref{eq:pmon}, we can conclude 
$$ \int_{B_r} |\gr u_j -\gr u|^p\dx \leq c\eps \|u\|_{W^{1,p}(N)}^p,\quad \forall\ j\geq j_0(\eps),$$
for any $B_r\subset \Om_j\cap N$. The continuity of the gradients $\gr u_j, \gr u$ from \eqref{eq:c1reg} implies that all points are Lebesgue points and hence the above is enough to conclude the following pointwise convergence of $\gr u_j -\gr u$; the estimate of \eqref{eq:c1reg} for $\gr u_j$ shows equicontinuity of 
the family $\{\gr u_j\}_{j\geq j_0}$ and uniform convergence $\gr u_j\to \gr u$ in $N$ as $j\to \infty$ follows from the Arzel\`a-Ascoli theorem. 
This used in \eqref{eq:wcd} directly leads us to conclude
$$ \lim_{j\to \infty}\int_{\del\Om_j} \eta |\gr u_j|^{p-1} d\h^{n-1} = \int_{\del\Om} \eta |\gr u|^{p-1} d\h^{n-1}$$
for any $\eta\in C^1(\bar N)$ and hence also for any $\eta\in C(\bar N)$ by taking approximations. The proof is finished. 
\end{proof}

\begin{remark}
The proof of weak convergence of $p$-harmonic measures above is easier compared to other measures with a higher exponent on the gradient because in our case, 
the exponent $(p-1)$ is precisely what appears in the weak form of the equations and so we do not require a reverse H\"older inequality as in \cite{CNSXYZ}. It is so remarked in \cite{J2} that the weak convergence in case of harmonic measure in \cite{J1} is easier than that of capacitary measures in \cite{J2} for $p=2$. 
\end{remark}

 Let $\mu_{\Om_j}= (\g_{\Om_j})_*\om_{p,j}$
be the corresponding measure as an approximation to $\mu_\Om$. As $\g_{\Om_j}\to \g_{\Om}$ $\h^{n-1}$-a.e. we have the weak convergence 
$\mu_{\Om_j}\wto \mu_{\Om}$ as a consequence of the above. Thus, 
\begin{equation}\label{eq:wmuconv}
 \lim_{j\to \infty}\int_{\Snn} w\, d\mu_{\Om_j}  = \int_{\Snn} w \, d\mu_{\Om}, \qquad \forall\ w\in C(\Snn).
\end{equation}


More generally, we consider $u=u_K\in W^{1,p}(N)$ for a general convex set $K$ and a neighborhood $N$ of $\del K$, where either 
$u$ satisfies \eqref{eq:omdir} if $K=\bar \Om$ or $u$ satisfies \eqref{eq:omdiremp} if $K=\del  K$. Then $\om_p$ defined by 
$d\om_p=|\gr u|^{p-1}d\h^{n-1}$ on $\del K$ and $\mu_K$ is as in \eqref{eq:mukint} and
defined with respect to \eqref{eq:gmap} and \eqref{eq:ginv} in general. We consider $\{K_j\}$ such that $d_\h(K_j,K)\to 0^+$. 
If $K_j=\bar\Om_j$ are of non-empty interior then we define $u_j\in W^{1,p}(\Om_j\cap N)$ similarly as \eqref{eq:omjdir},  
otherwise if $K_j=\del K_j$ and $K=\del K$ then we define $u_j\in W^{1,p}(N)$ as the unique solution of 
\begin{equation}\label{eq:omjdiremp}
 \begin{cases}
  \dv \big(|\gr u_j|^{p-2}\gr u_j\big)=0,\ \ &\text{in}\ N\mns K_j;\\
  u_j=0, \ \  &\text{in}\ K_j\cap  N;\\
  u_j= u,\ \  & \text{in}\ \del N\mns K_j.  
 \end{cases}
\end{equation}
The $p$-harmonic measures are defined by $d\om_{p,j}=|\gr u_j|^{p-1}d\h^{n-1}\on_{\del K_j}$ and the weak convergence 
$\om_{p,j}\wto \om_p$ can be obtained similarly as Proposition \ref{prop:pconvmsr}. 
Therefore, the weak convergence for measures for the $p$-harmonic Minkowski problem on general convex sets,  
\begin{equation}\label{eq:wconph}
\mu_{K_j} \wto \mu_K, \qquad\qquad \text{if}\quad d_\h(K_j,K)\to 0^+,
\end{equation}
when the measures are defined as in \eqref{eq:mukint} with respect to \eqref{eq:gmap} and \eqref{eq:ginv}, 
can also be obtained using the above and approximation of general convex sets by polytopes, 
in ways similar to that of \eqref{eq:wcons}. We refer to \cite[Chapter 4,5]{Sc} for details of such arguments.  

In the following, given a convex set $K$, we denote the functional 
\begin{equation}\label{eq:Gfunc}
\Gamma (K):= \int_{\Snn} h_K ( \xi)\, d \mu_K( \xi) 
\end{equation}
and $\Gamma(\Om)=\Gamma(K)$ if $\bar \Om=K$. It is easy to see from \eqref{eq:hdh} and \eqref{eq:wconph}, that $\Gamma(K_j)\to \Gamma(K)$ uniformly if 
$d_\h(K_j, K)\to 0^+$ as $j\to \infty$. 

\begin{remark}
For a compact convex set $K$ of non-empty interior note that 
$$\Gamma (K)= \int_{\Snn} h_K ( \xi)| \gr u(\inv{\g_K}(\xi))|^{p-1} d S_K( \xi) $$ closely resembles the well known formula 
of volumes given by 
$$ |K|=\frac{1}{n} \int_{\Snn} h_K ( \xi) d S_K( \xi).$$
Hence we can conclude that if $K\mapsto \mu_K$ is defined by a common $p$-harmonic function 
$u$ for $C^2_+$ domains then $\Gamma$ is a measure absolutely continuous with respect to the Lebesgue measure. For 
non-smooth domains or convex sets of empty interior, the measure $\Gamma$ can be singular. 
\end{remark}

We complete this subsection with the following. 
\begin{lemma}\label{prop:gomtconv}
Given convex domains $\Om_j, \Om_{j,0}\subset \R^n$ with
$d_\h(\Om_j,\Om), d_\h(\Om_{j,0},\Om_0) \to 0^+$ as $j\to \infty$, 
let $\Om^t=\Om+t\Om_0$ and $\Om_j^t=\Om_j+t\Om_{j,0}$ for $|t|<\tau$ small enough. 
Then, let $$d\mu_{\Om^t}= |\gr u (\cdot, t)|^{p-1}d\h^{n-1}\on_{\del \Om^t},$$ where 
$u(\cdot, t)$ is the unique solution of the Dirichlet problem \eqref{eq:omtdir} and let $d\mu_{\Om_j^t}$ be defined similarly; then with $\Gamma$
as \eqref{eq:Gfunc}, if $t\mapsto \Gamma (\Om_j^t)$ is differentiable at $t=0$ for all $j\geq j_0$ then 
$t\mapsto \Gamma (\Om^t)$ is also differentiable at $t=0$ and we have
\begin{equation}\label{eq:gomtconv}
\lim_{j\to \infty}\frac{d}{dt}\Big|_{t=0} \Gamma (\Om_j^t) = \frac{d}{dt}\Big|_{t=0} \Gamma (\Om^t)= \frac{d}{dt}\Big|_{t=0} \lim_{j\to \infty} \Gamma (\Om_j^t).
\end{equation}
\end{lemma}
\begin{proof}
From \eqref{eq:suppminkadd}, we note that 
\begin{align*}
\Gamma (\Om_j^t)-\Gamma (\Om_j) &= \int_{\Snn} h_{\Om_j^t} \, d \mu_{\Om_j^t} - \int_{\Snn} h_{\Om_j} \, d \mu_{\Om_j}\\
&=  \int_{\Snn} h_{\Om_j} d \mu_{\Om_j^t}-\int_{\Snn} h_{\Om_j} \, d \mu_{\Om_j} + t\int_{\Snn}  h_{\Om_{j,0}} \, d \mu_{\Om_j^t}. 
\end{align*}
Recalling \eqref{eq:hdh}, we know that $h_{\Om_j},  h_{\Om_{j,0}} \to  h_{\Om},  h_{\Om_{0}}$ uniformly as $j\to \infty$. 
Hence for each $t\neq 0$, we have $(\Gamma (\Om_j^t)-\Gamma (\Om_j))/t \to (\Gamma (\Om^t)-\Gamma (\Om))/t$ uniformly (independent of $t$) as $j\to \infty$.
Also from \eqref{eq:hdh}, note that $$d_\h (\Om_j^t, \Om_j)= |t|\| h_{\Om_{j,0}}\|_{L^\infty}\qquad\text{and}
\qquad d_\h (\Om^t, \Om)= |t|\| h_{\Om_{0}}\|_{L^\infty},$$ which converge to 
zero uniformly (independent of $j$) as $t\to 0$ as well. From \eqref{eq:wmuconv} and above, $\mu_{\Om_j^t}\wto \mu_{\Om^t}$ as $j\to \infty$ for all $t$ and 
$\mu_{\Om_j^t}\wto \mu_{\Om_j}$ as $t\to 0$ for all $j\geq j_0$. 
This is enough to conclude the proof. 
\end{proof}

\subsection{Proof of existence} We shall prove the existence of $p$-harmonic measures prescribing measures on spheres in two steps; first, for discrete measures and then 
for the general case thereby leading to the proof of Theorem \ref{thm:main} in the end. 

The following proposition resembles the Hadamard's variational formula. 
\begin{proposition}\label{prop:gomt}
Given convex domains $\Om^t=\Om+t\Om_0$ with $|t|<\tau$ small enough, $v= h_{\Om_0}>0$ and $\Gamma$ as in \eqref{eq:Gfunc} with 
$d\mu_{\Om^t}= |\gr u (\cdot, t)|^{p-1}d\h^{n-1}\on_{\del \Om^t}$ where 
$u(\cdot, t)$ is the unique solution of the Dirichlet problem \eqref{eq:omtdir}, we have 
\begin{equation}\label{eq:gomt}
\frac{d}{dt}\Big|_{t=0} \Gamma (\Om^t)= (n-p+1) \int_{\Snn} v \,d\mu_{\Om}.
\end{equation}
\end{proposition}
\begin{proof}
First, we assume that $\Om, \Om_0$ are of class $C^{2,\alpha}_+$. Then, recalling $d\mu_{\Om^t} ( \xi)= \F[h+tv]( \xi)\,d \xi$ where $h_\Om= h$, 
observe that 
 \begin{equation}\label{constrgam}
\begin{aligned}
\frac{d}{dt}\Big|_{t=0} \Gamma (\Om^t)&= \frac{d}{dt}\Big|_{t=0} \int_{\Snn} h_{\Om^t} ( \xi) \,d\mu_{\Om^t}( \xi)
= \frac{d}{dt}\Big|_{t=0} \int_{\Snn} (h+tv) \F[h+tv]( \xi) \,d  \xi\\
&= \int_{\Snn} v( \xi) \F[h]( \xi) \,d  \xi + \int_{\Snn} h( \xi) \frac{d}{dt}\Big|_{t=0}  \F[h+tv]( \xi) \,d  \xi\\
&= \int_{\Snn} v( \xi) \,d\mu_{\Om} ( \xi) + \int_{\Snn} h( \xi) \LL_{h}[v]( \xi)\,d  \xi\\
&= \int_{\Snn} v( \xi) \,d\mu_{\Om} ( \xi) + \int_{\Snn} v( \xi) \LL_{h}[h]( \xi)\,d  \xi,
\end{aligned}
\end{equation}
where the self-adjointness of 
Corollary \ref{cor:llvsym} is used to obtain the last equality of the above. Now, recalling \eqref{eq:flhh} of Lemma \ref{lem:flhh}, we have  
$\LL_{h}[h]( \xi)=(n-p)\F[h]( \xi)$.  Using this on \eqref{constrgam}, we obtain
$$ \frac{d}{dt}\Big|_{t=0} \Gamma (\Om^t) = (n-p+1) \int_{\Snn} v( \xi) \,d\mu_{\Om} ( \xi),$$
for the case of $C^{2,\alpha}_+$ domains. The general case follows from approximation of convex domains by $C^{2,\alpha}_+$ domains using 
Theorem \ref{thm:c2approx} and \eqref{eq:wconvmsr} along with \eqref{eq:wmuconv} and \eqref{eq:gomtconv}. 
\end{proof}

In the following lemma, we consider the possibility of having a lower dimensional convex set
in order to be able to assert the existence theorem for convex sets of 
non-empty interior. 
\begin{lemma}\label{lem:lowdim}
Given a compact convex set $K$ and a measure $\mu$ on $\Snn$, we have the following: 
\begin{enumerate}
\item If $\dim_\h (K)= n-1$, then $K$ is contained in a hyperplane; precisely, there exists $\xi_0\in \Snn$ and $c\in \R$ such that $K\subset \{x: \inp{x}{\xi_0}= c\}$ and $\mu_K$ as in \eqref{eq:mukint} is given by 
$$ \mu_K = \Big(\int_{K} |\gr u|^{p-1} d\h^{n-1}\Big) \big(\delta_{\xi_0} +\delta_{-\xi_0}\big). $$
\item If $\mu= a\big(\delta_{\xi_0} +\delta_{-\xi_0}\big)$ for some $\xi_0\in \Snn$ and $a>0$, then $\mu$ does not satisfy $(i)$ of \eqref{eq:excond}.  
\end{enumerate}
\end{lemma}

\begin{proof}
To prove $(1)$, recalling \eqref{eq:kinthk} since $K=\cap_{\xi\in \Snn} H_{\xi, h_K(\xi)}^-$, observe that $\dim_\h (K)= n-1$ can only occur if at least two half-spaces share 
a common boundary and $K$ is degenerated to a $1$-codimensional subset of a hyperplane, i.e. $H_K(\xi)=H_K(\xi')$ for some 
$\xi\neq\xi'\in \Snn$ which can happen only if $\xi'=-\xi $ and 
$h_K(\xi')=-h_K(\xi) $. These together imply $h_K(-\xi)=-h_K(\xi)$; but, from 
the definition of support function, since 
$$-h_K (-\xi)\leq \inp{x}{\xi}\leq h_K (\xi), \qquad \forall\ x\in K, \xi\in \Snn,$$ therefore equality has to hold. 
Thus, we have $K\subset \{x: \inp{x}{\xi_0}= c\}$ where $\xi_0=\xi$ and $ c= h_K(\xi)$. Then, we observe that, in the sense of 
\eqref{eq:gmap}, we have $\g_K(x)=\{\xi_0,-\xi_0\}$ for $\h^{n-1}$-a.e. $x\in K$. 
Precisely, 
if $F\subset K$ is the lower dimensional boundary of $K$ within the hyperplane then $\g_K(x)=\{\xi_0,-\xi_0\}$ for all $x\in K\mns F$ and $\h^{n-1}(F)=0$. 
Also in this case $K=\del K$. 
Hence, for any measurable $E\sub \Snn$, 
note that, $\inv{g_K}(E)$ in the sense of \eqref{eq:ginv} is of full $\h^{n-1}$-measure only if 
$\xi_0\in E$ or $-\xi_0\in E$; to be precise, letting $F'=\set{x\in F}{\g_K(x)\cap E\neq \emp}\sub F$, we have that 
$\inv{g_K}(E)= (K\mns F) \cup F'$ if $\{\xi_0,-\xi_0\}\cap E\neq \emp$ and otherwise $\inv{g_K}(E)= F'$. 
Thus, 
$$ S_K= \h^{n-1}(K) \big(\delta_{\xi_0} +\delta_{-\xi_0}\big) \quad\text{and}\quad  \mu_K = \Big(\int_{K} |\gr u|^{p-1} d\h^{n-1}\Big) \big(\delta_{\xi_0} +\delta_{-\xi_0}\big),$$
from \eqref{eq:mukint}, which completes the first part. 

To prove $(2)$, we observe that if $\mu= a\big(\delta_{\xi_0} +\delta_{-\xi_0}\big)$ for some $\xi_0\in \Snn$ and $a>0$, then 
$$\int_{ \mathbb{S}^{n-1}}  | \inp{\zeta}{ \xi} | \, d \mu (  \xi ) 
= a \Big( | \inp{\zeta}{ \xi_0} |+ |\inp{\zeta}{ -\xi_0} |\Big)=2a | \inp{\zeta}{ \xi_0} |, $$
which vanishes for any $\zeta\in \Snn$ normal to $\xi_0$. This completes the second part and the proof. 
\end{proof}

The following theorem is a discrete version of Theorem \ref{thm:main}, which is an essential step towards the general existence theorem. Here we require a technical assumption of an antipodal condition and $p\neq n+1$, which shall be ultimately removed. 

\begin{theorem}\label{thm:exdiscrete}
Let $ \mu $ be a discrete positive Borel measure on $\mathbb{S}^{n-1}$ satisfying conditions \eqref{eq:excond} and has no antipodal pair of point masses i.e. if $ \mu ( \{ \xi\})>0$ for some  $  \xi \in \mathbb{S}^{n-1}$ then $ \mu (\{- \xi\}) =0$. 
Then there exists a polytope $P$ with non-empty interior such that $\mu_{P} = \mu$ and $\Gamma(P)=1$ with $\mu_{P}= (\g_P)_* \om_p$ as in 
\eqref{eq:mukint}, where $\om_p$ is any $p$-harmonic measure on $\del P$ and $\Gamma$ is as in \eqref{eq:Gfunc}, for any $1<p<\infty$ with $p\neq n+1$. 
\end{theorem}

\begin{proof}
By virtue of the duality relation of Theorem \ref{thm:possupp} between convex domains and positive continuous functions on $\Snn$, 
the existence of the polytope $P$ shall be obtained from 
existence of minimizer of the following constrained minimization problem 
\begin{equation}\label{eq:minprobh}
b= \inf\set{\int_{\Snn} h \, d\mu }{h\in C_+(\Snn),\, \Gamma(\Om_h)\geq 1}, 
\end{equation}
where $\Om_h$ is the Wulff shape of $h$ as in \eqref{eq:wshape}, hence 
$\bar\Om_h= \{x\in \R^n : \inp{x}{ \xi}\leq  h( \xi)\ \forall\  \xi\in \Snn\}$, and $\Gamma$ is as in \eqref{eq:Gfunc} for a suitably defined measure $\mu_{\Om_h}$, to be mentioned below. 

Since $\mu $ is discrete, there exists $ \xi_1,\ldots,  \xi_m\in \Snn$ such that we have 
$$\mu= \sum_{i=1}^m c_i \delta_{ \xi_i}$$ where, the condition 
\eqref{eq:excond} implies for all $ \xi\in \Snn$, we have
\begin{equation}\label{eq:excondiscr}
\sum_{i=1}^m c_i |\inp{ \xi}{ \xi_i}|>0\quad\text{and}\quad \sum_{i=1}^m c_i \xi_i =0, 
\end{equation}
and in particular, $c_1,\dots,c_m >0$. 
Let us define $P_1 = \big\{x\in \R^n : \inp{x}{ \xi_i}\leq h( \xi_i)\ \forall\ i\in\{1,\ldots, m\}\big\}$ for some $h\in C_+(\Snn)$,  
which is a polytope from Theorem \ref{thm:charpoly}. We note that $\bar\Om_h\sub P_1$ 
and 
$h_{P_1}( \xi_i) = \sup_{x\in {P_1}} \inp{x}{ \xi_i} \leq h( \xi_i)$ from the definition of $P_1$ and $\h^{n-1}(\del\Om_h)\leq \h^{n-1}(\del P_1)$. Thus, we have 
\begin{equation}\label{p1e}
\int_{\Snn} h_{P_1} \,d\mu = \sum_{i=1}^m c_i h_{P_1}( \xi_i) \leq \sum_{i=1}^m c_i h( \xi_i) = \int_{\Snn} h \,d\mu .
\end{equation}
We also note that $0\in \Om_h$ if $h> 0$ and from Theorem \ref{thm:possupp}, for $S_{\Om_h}$-a.e. $\xi\in \Snn$, 
we have $h(\xi) =h_{\Om_h}(\xi)\leq h_{P_1}(\xi)$ that may or may not include the set $\{\xi_1,\ldots,\xi_m\}$. 

For any bounded sequence $h_j\in C_+(\Snn)$, let us denote the corresponding sequence of polytopes as above, i.e. 
$P_j = \big\{x\in \R^n : \inp{x}{ \xi_i}\leq h_j( \xi_i)\ \forall\ i\in\{1,\ldots, m\}\big\}$. 
Then, $\bar\Om_{h_j}\sub P_j$ and 
\eqref{p1e} leads to 
$$ \sum_{i=1}^m c_i h_{P_j}( \xi_i) \leq \limsup_{j\to \infty} \int_{\Snn} h_j \,d\mu =: M<\infty $$ 
and since $c_i>0$, this implies for all $x\in P_j$ we have $\inp{x}{ \xi_i}\leq h_{P_j}( \xi_i)\leq M/(\min_{1\leq i\leq m} c_i)$. 
Thus, we can conclude that the sequence 
of polytopes $\{P_j\}$ is bounded 
and hence, from Blaschke selection theorem (Theorem \ref{thm:bla}), 
there exists $P$ such that, up to subsequence, $d_\h (P_j, P)\to 0^+$ as $j\to \infty$. 
From Theorem \ref{thm:charpoly}, $P$ is also a polytope with normals $\xi_1,\ldots, \xi_m$. 

Let $u=u_P$ be the $p$-harmonic function corresponding to any fixed $p$-harmonic measure of $P$, 
defined in the neighborhood of $\del P$ that contains $\del P_j$ for all $j\geq j_0$, since from \eqref{eq:hdh} 
$d_\h (\del P_j, \del P)\to 0^+$ as $j\to \infty$. With respect to this $u$, we choose $u_j$ uniquely as solutions 
of the Dirichlet problem \eqref{eq:omjdir} as in the previous subsection, so that we have $\mu_{P_j}\wto \mu_{P}$ from \eqref{eq:wconph}.
We shall define $\mu_{\Om_{h_j}}$ with respect to $u_j$ similarly as $\mu_{P_j}$ on neighborhoods of 
$\del P_j$ that contain $\del\Om_{h_j}$ for $j\geq j_0$ large enough. In other words, we have defined 
$$\mu_{P_j}=(\g_{ P_j})_*\big( |\gr u_j|^{p-1}\h^{n-1}\on_{\del P_j}\big)\quad\text{and}\quad 
\mu_{\Om_{h_j}}= (\g_{\Om_{h_j}})_*\big(  |\gr u_j|^{p-1}\h^{n-1}\on_{\del \Om_{h_j}}\big).$$
Now, let $h_j\in C_+(\Snn)$ be a minimizing sequence of the infima \eqref{eq:minprobh} where $\Gamma$ is as in \eqref{eq:Gfunc} with
$\mu_{\Om_{h_j}}$ as above. 
Since $\mu_{P_j}\wto \mu_{P}$, we have 
hence $\Gamma (P_j)\to \Gamma (P)$ as $j\to \infty$ and  
since $\Om_{h_j}\sub P_j$, hence from \eqref{eq:minprobh} note that $1\leq \Gamma (\Om_{h_j})\leq \Gamma (P_j)$ and therefore, 
$$\int_{\Snn} h_P d\mu_P=\Gamma (P) = \lim_{j\to \infty} \Gamma (P_{j}) \geq 1.$$

Suppose, $P$ is of empty interior, hence $P=\del P$ and Hausdorff dimension $\dim_\h (P)<n$. 
Notice that the above enforces $\dim_\h (P)\geq n-1$ (otherwise, $\h^{n-1}(P)=0=\h^{n-1}(\del P)$ which, from \eqref{eq:mukint}, implies $\mu_P\equiv 0$ and hence
$\Gamma(P)=0$ contradicting the above). 
Then, $\dim_\h(P)=n-1$.  
As in $(1)$ of Lemma \ref{lem:lowdim}, this can occur only if at least two half-spaces of $P$ share 
a common boundary and degenerate the whole polytope to a $1$-codimensional face; if those are the half spaces corresponding to $\xi_i\neq\xi_j\in \Snn$, it enforces $\xi_j= -\xi_i$. But since $c_i, c_j>0$, this contradicts the antipodal condition. Thus, $P$ is of non-empty interior. 

Therefore, let $P_0\neq \emp$ be the interior of $P$ and hence $\Gamma(P_0)\geq 1$, $h_{P_0}$ is a minimizer of \eqref{eq:minprobh} and from $(ii)$ of \eqref{eq:excond}, the function $ \xi \mapsto h_{P_0}( \xi)+ \gamma  \xi$ is also a minimizer
for any $\gamma>0$. Hence, without loss of generality, we can regard $h_{P_0}( \xi)>0$ for all $ \xi\in \Snn$. 
Recalling Theorem \ref{thm:possupp}, 
for any arbitrary $v\in C_+(\Snn)$, let $\Om_v$ is the Wulff shape of $v$ as in \eqref{eq:wshape} and then $h_{\Om_v}=v$ holds 
$S_{\Om_v}$-a.e. 
Let $h_0= h_{P_0}$ and 
$$\Om^t= P_0+t\Om_v \quad \text{so that}\quad h_{\Om^t}= h_{P_0} +t h_{\Om_v}=h_0+tv,\ \ S_{\Om_v}-a.e.$$ 
With respect to $u=u_P$ denoted above that is $p$-harmonic in a neighborhood of $\del P=\del P_0$, we choose 
$\mu_{\Om^t}$ to denote $\Gamma (\Om^t)$ as in \eqref{eq:Gfunc}, defined by 
$d\mu_{\Om^t}= |\gr u (\cdot, t)|^{p-1}d\h^{n-1}\on_{\del \Om^t}$ where 
$u(\cdot, t)$ is the unique solution of the Dirichlet problem \eqref{eq:omtdir}. 
Without loss of generality, we can replace $P_j$ and $\Om_{h_j}$ above with infinitesimal dilations by $(1+ t)$ with $|t|<1/j$ which does not affect the Hausdorff convergence. 
Then, $\mu_{\Om^t}$ with respect to $u (\cdot, t)$ for $|t|<1/j$ forms a special class of approximants among $\mu_{\Om_{h_j}}$ with respect to $u_j$ of the above (see Remark \ref{rem:conn}). Therefore, the minimization at $h_0$ leads to the existence of Lagrange multiplier $\lambda$ such that we have
\begin{equation}\label{constr}
 \frac{d}{dt}\Big|_{t=0} \int_{\Snn} (h_0+tv) \, d\mu = \lambda \frac{d}{dt}\Big|_{t=0} \Gamma (\Om^t). 
 \end{equation}
 Now \eqref{eq:gomt} of Proposition \ref{prop:gomt} together with \eqref{constr}, lead to 
$$ \int_{\Snn} v\,d\mu = \lambda(n-p+1)\int_{\Snn} v\,d\mu_{P_0}$$
for any $v\in C_+(\Snn)$. Then, the above is also obtained for any $v\in C(\Snn)$ using $v= v^+ - v^-$, which is enough to conclude
$ \mu= \lambda(n-p+1)\mu_{P_0}$. Since $h_{0}$ is the minimizer of \eqref{eq:minprobh}, we have 
$$ b= \int_{\Snn} h_{0} d\mu= \lambda(n-p+1)\int_{\Snn} h_{P_0} d\mu_{P_0} =\lambda(n-p+1) \Gamma(P_0),$$
so that $\lambda = b/(n-p+1) \Gamma(P_0)$ and hence, $\mu= \frac{b}{\Gamma(P_0)}\mu_{P_0}$. In other words, we have  
$$ \frac{\mu}{\displaystyle \int_{\Snn} h_{P_0} d\mu}= \frac{\mu_{P_0}}{\displaystyle \int_{\Snn} h_{P_0} d\mu_{P_0}},$$
which is equivalent to $\mu=c\mu_{P_0}$ for some $c>0$. 
By rescaling $P_0$, we can ensure $\mu_{P_0}= \mu$ and $\Gamma(P_0)=1$. The proof is complete. 
\end{proof}

Finally, we are ready to prove the existence theorem for general Borel measures. 
\begin{proof}[Proof of Theorem \ref{thm:main}]
Given the measure $\mu$ on $\Snn$ satisfying \eqref{eq:excond}, note that from $(i)$ of \eqref{eq:excond}, we have 
$$ \delta=\inf_{\zeta\in \Snn} \int_{ \mathbb{S}^{n-1}}  | \inp{\zeta}{ \xi} | \, d \mu (  \xi )  >  0. $$
There exists discrete measures $\mu_j$ on $\Snn$ such that $\mu_j\wto \mu$ as $j\to \infty$ and we can assume that there exists $c=c(\delta)>0$ such that 
for all $j\in \N$, we have 
\begin{equation}\label{eq:jupbd}
\inf_{\zeta\in \Snn} \int_{ \mathbb{S}^{n-1}}  | \inp{\zeta}{ \xi} | \, d \mu_j (  \xi )  \geq c.
\end{equation}

Without loss of generality we can assume that $\mu_j$'s have no antipodal pair of point masses 
Indeed, if $ \mu_j ( \{ \xi_j\}), \mu_j ( \{ -\xi_j\})>0$ then $\mu_j= c(\delta_{\xi_j} +\delta_{-\xi_j}) +\nu_j$ for some other discrete measure $\nu_j$ and we can replace $\mu_j$ with 
$\mu_j^{\eps_j}=  c(\delta_{\xi'_j} +\delta_{-\xi_j}) +\nu_j$ where $0<|\xi'_j-\xi_j|<\eps_j$ for $\eps_j>0$ small enough, so that $\mu_j^{\eps_j} ( \{ \xi_j\})=0$ while $\mu_j^{\eps_j} ( \{ -\xi_j\})>0$ and 
$\mu_j^{\eps_j}\wto \mu$ as $j\to \infty$. After finitely many replacements by such small perturbations on all antipodal points, we can regard $\mu_j$'s have no antipodal point masses for each $j\in \N$. 

Now we can use Theorem \ref{thm:exdiscrete} above to conclude for all $j\in \N$ there exists a polytope $P_j$ of non-empty interior such that $\mu_{P_j}= \mu_j$ and 
$\Gamma(P_j)=1$ defined as in \eqref{eq:Gfunc} with respect to $\mu_{P_j}= (\g_{P_j})_* \om_{p,j}$ where $\om_{p,j}$ is any $p$-harmonic measure on $\del P_j$ for $p\neq n+1$. 
For the case $p=n+1$, we denote $\mu_{P_j}^{\eps_j}$ and $\Gamma^{\eps_j}$ corresponding to $p=n+1+\eps_j$ and using Theorem \ref{thm:exdiscrete}, take $P_j$ such that $\mu_{P_j}^{\eps_j}= \mu_j$ and 
$\Gamma^{\eps_j}(P_j)=1$
where $\eps_j\neq 0$ is small enough; the limits remain the same since $\|\mu_{P_j}^{\eps_j}-\mu_{P_j}\|,\ |\Gamma^{\eps_j}(\cdot)-\Gamma (\cdot)|\to 0^+$ 
as $\eps_j\to 0$. 
We show that the sequence of polytopes $\{P_j\}$ is bounded. Indeed, for any $x\in \bar P_j$ we have 
$$ \int_{\Snn} \inp{x}{\xi} d \mu_j (\xi) \leq \int_{\Snn} h_{P_j}(\xi) d \mu_j (\xi) 
= \int_{\Snn} h_{P_j} d\mu_{P_j}=\Gamma(P_j)=1.$$
Up to a suitable translation, we have $0\in P_j$ and $x_j, -x_j\in \del P_j$ such that $\diam(P_j)=2|x_j|$. Therefore, we have 
$\langle \pm x_j, \xi\rangle\leq h_{P_j}(\xi)$ and hence $|\inp{x_j}{\xi}|\leq h_{P_j}(\xi)$. Using this on the above 
 together with \eqref{eq:jupbd}, we get
 $$ c|x_j|\leq \int_{\Snn} |\inp{x_j}{\xi}| d \mu_j (\xi) \leq \int_{\Snn} h_{P_j}(\xi) d \mu_j (\xi) =1,$$
which implies $\diam(P_j)/2=|x_j|\leq 1/c$. From Blaschke selection theorem (Theorem \ref{thm:bla}), 
there exists a convex set $K$ such that, up to subsequence, $d_\h (P_j, K)\to 0^+$ as $j\to \infty$.
Let $\om_p$ be any $p$-harmonic measure on $\del K$ with respect to which we have $\mu_K$ as in \eqref{eq:mukint} and $u=u_K$ be the corresponding 
$p$-harmonic function in a neighborhood of $\del K$. Then, we choose the $p$-harmonic measures $\om_{p,j}$ as 
$d\om_{p,j}= |\gr u_j|^{p-1}d\h^{n-1}\on_{\del P_j}$ 
where, with respect to this $u=u_K$, we choose $u_j$'s uniquely as the solutions 
of the Dirichlet problem \eqref{eq:omjdir} as in the previous subsection, so that we have $\mu_{P_j}\wto \mu_{K}$ from \eqref{eq:wconph} and 
hence $\Gamma (P_j)\to \Gamma (K)$ as $j\to \infty$, leading to $\Gamma(K)=1$ and from uniqueness of the weak limit, $\mu_K=\mu$. 
Now we show that $K$ has a non-empty interior. 

Suppose $K$ is of empty interior so that $K=\del K$ then, the Hausdorff dimension $\dim_\h(K)<n$. Now if $\dim_\h(K)< n-1$ then $\h^{n-1}(K)=0=\h^{n-1}(\del K)$ would imply $\mu_K\equiv 0$ from \eqref{eq:mukint} and hence $\Gamma(K)=0$ which contradicts $\Gamma(K)=1$. Therefore, as $\dim_\h(K)\in \N$, we conclude $\dim_\h(K)= n-1$.
Since $\mu_K=\mu$ satisfy \eqref{eq:excond}, we invoke Lemma \ref{lem:lowdim} for a contradiction. 

Therefore, there exists a domain $\Om$ such that $K=\bar\Om$ and $\mu_\Om= \mu$. The proof is complete. 
\end{proof}

\addtocontents{toc}{\SkipTocEntry}
\section*{Acknowledgement} The authors are thankful to Peter Van Hintum and the anonymous referee for their helpful remarks on a previous version of this paper. 

\addtocontents{toc}{\SkipTocEntry}
\section*{Data availability}
No data was used for the research described in the article.

\bibliographystyle{plain}
\bibliography{minkprob}

\end{document}